\documentclass [12pt] {report}
\usepackage[dvips]{graphicx}
\usepackage{epsfig, latexsym,graphicx}
\usepackage{amssymb}
\newcommand {\eqdef} {\stackrel{\rm def}{=}}
\newcommand {\D}[2] {\displaystyle\frac{\partial{#1}}{\partial{#2}}}

\newcommand {\Dd}[3] {\displaystyle\frac{\partial^2{#1}}{\partial{#2}\partial{#3}}}
\newcommand {\al} {\alpha}

\newcommand {\ga} {\gamma}
\newcommand {\la} {\lambda}

\newcommand {\si} {\sigma}

\newcommand {\de} {\delta}
\newcommand {\prtl} {\partial}
\newcommand {\fr} {\displaystyle\frac}
\newcommand {\wt} {\widetilde}
\newcommand {\be} {\begin{equation}}
\newcommand {\ee} {\end{equation}}
\newcommand {\ba} {\begin{array}}
\newcommand {\ea} {\end{array}}
\newcommand {\bp} {\begin{picture}}
\newcommand {\ep} {\end{picture}}
\newcommand {\bc} {\begin{center}}
\newcommand {\ec} {\end{center}}
\newcommand {\bt} {\begin{tabular}}
\newcommand {\et} {\end{tabular}}
\newcommand {\lf} {\left}
\newcommand {\rg} {\right}

\newcommand {\cC} {{\cal C}}
\newcommand {\cF} {{\cal F}}

\newcommand {\cR} {{\cal R}}
\newcommand {\cS} {{\cal S}}

\newcommand {\ses} {\medskip}

\newcommand {\e} {\mathop{\rm e}\nolimits}

\newcommand {\bibit} {\bibitem}
\newcommand {\nin} {\noindent}

\newcommand {\Rho} {\mbox{\large$\rho$}}

\usepackage{epsfig, latexsym,graphicx}
 \usepackage{amsmath}

\newcommand {\cD} {{\cal D}}

\def\2#1#2#3{{#1}_{#2}\hspace{0pt}^{#3}}
\def\3#1#2#3#4{{#1}_{#2}\hspace{0pt}^{#3}\hspace{0pt}_{#4}}
\newcounter{sctn}
\def\sec#1.#2\par{\setcounter{sctn}{#1}\setcounter{equation}{0}
                  \noindent{\bf\boldmath#1.#2}\bigskip\par}

\begin {document}

\begin {titlepage}

\vspace{0.1in}

\begin{center}

{\large \bf   Finsler connection preserving angle in dimensions $N\ge3$}

\end{center}

\vspace{0.3in}

\begin{center}

\vspace{.15in} {\large G.S. Asanov\\} \vspace{.25in}
{\it Division of Theoretical Physics, Moscow State University\\
119992 Moscow, Russia\\
{\rm (}e-mail: asanov@newmail.ru{\rm )}} \vspace{.05in}

\end{center}

\begin{abstract}

 We show that if a
Finsler space
is  conformally automorphic to a Riemannian space
and the automorphism is positively homogeneous with respect to tangent vectors,
then  the indicatrix of the Finsler space is a space of constant curvature.
In this case, the  Finslerian two-vector angle can explicitly be found,
which gives rise to simple and explicit representation for
the  connection preserving the  angle
in the indicatrix-homogeneous case.
The  connection is metrical
and
the  Finsler space is obtainable from
the Riemannian space
by means of the parallel deformation.
Since also the transitivity of covariant derivative holds,
in such Finsler spaces
 the metrical non-linear angle-preserving connection
 is  the respective
 export of
 the
  metrical  linear Riemannian  connection.
From the commutators of covariant derivatives the associated
curvature tensor is found.
In case of the
${\cal FS}$-space,
the  explicit  example
  of the conformally automorphic transformation
can be  developed, which entails the explicit connection coefficients
and the metric function of the Finsleroid type.

\end{abstract}

\end{titlepage}

\vskip 1cm

\ses

\ses

\setcounter{sctn}{1} \setcounter{equation}{0}

\bc  {\bf 1. Motivation and description}  \ec

\ses

\ses

In any dimension $N\ge3$
the Finsler metric function $F$
geometrizes the tangent bundle $TM$ over the base manifold $M$
such that
at each point $x\in M$
the tangent space $T_xM$ is  endowed with the curvature tensor
constructed from the respective Finslerian metric tensor  $g_{\{x\}}(y)$
by means of the conventional rule of the Riemannian geometry
considering $y$ to be the variable argument.
There arises the Riemannian space
$\cR_{\{x\}}=\{g_{\{x\}}(y),  T_xM\}$ supported by the point $x\in M$
 such that $T_xM$ plays the role of the base manifold for the space.
In the Riemannian limit of the  Finsler space,
 the spaces
 $\cR_{\{x\}}$
 are Euclidean spaces,
so that
the tensor
 $g_{\{x\}}(y)$ is independent of $y$.
 The conformally flat structure of the spaces
 $\cR_{\{x\}}$
can naturally  be taken to treat as  the next  level of generality
of the Finsler space.
Can the metrical connection preserving the two-vector angle be introduced
on that level?

The  deformation of the Riemannian space to the Finsler space
proves to be
 convenient invention to apply.
Namely, in the particular case when the Riemannian space
can be  deformated to  the Finsler space
characterized by
the  conformally flat structure of the spaces
 $\cR_{\{x\}}$
the positive and clear answer to the above question can be arrived at.

Given an $N$-dimensional  Riemannian space
${\cal R}^N=(M,S)$, where $S$ denotes the Riemannian metric function,
one may endeavor to  obtain a Finsler space
${\cal F}^N=(M, F)$ by applying an appropriate
{\it deformation}
 ${\cal C}$ of the space ${\cal R}^N$.
The notation $ F$  stands for the Finsler metric
function.
The  base manifold $M$ is keeping the same for both the spaces,
${\cal R}^N$ and ${\cal F}^N$.

We assume that the transformation ${\cal C}$ is {\it restrictive}, in the sense that no point $x\in M$ is shifted under the transformation, so that in each tangent space
$T_xM$ the deformation maps tangent vectors $y\in T_xM$ into the tangent vectors
of the same
$T_xM$:
\be
 y={\cal C}(x,\bar y), \quad y,\bar y\in T_xM.
\ee
In general, this transformation is non-linear with respect to $\bar y$. Non-singularity and
sufficient smoothness are always  implied.

We may evidence in the
Riemannian
space ${\cal R}^N$
the {\it metrical  linear Riemannian  connection}  ${\cal RL}$, which
in terms of local coordinates $\{x^i\}$ introduced in $M$
is given by
 \be
{\cal RL}=\{L^m{}_j,L^m{}_{ij}\}: \qquad    L^m{}_j=-a^m{}_{ij}y^i,
   \quad
   L^m{}_{ij}=a^m{}_{ij},
\ee
 with $a^m{}_{ij}=a^m{}_{ij}(x)$ standing for the Christoffel symbols
constructed from the Riemannian metric tensor $a_{mn}(x)$
of the space ${\cal R}^N$. The indices
$i,j,...$ are specified on the range $(1,\dots, N)$.
The respective covariant derivative $\nabla$ can be introduced
in the natural way, namely by means of the definition (4.18)
which uses the operator
\be
d^{\text{Riem}}_i=
\D{}{x^i}+L^k{}_i\D{}{y^k},
\ee
 to act on tensors considered on the tangent bundle
underlined the space ${\cal R}^N$.
In the space, the scalar product
${\langle y_1,y_2\rangle}^{\text{Riem}}_{\{x\}}=a_{mn}(x)y_1^my_2^n$
of two vectors $y_1,y_2$ supported by a fixed point $x\in M$
is linear with respect to each vector, which gives rise to  the profound meaning of the connection (1.2) to preserve the product
 under the entailed parallel transports
of the entered vectors along curves running on
 $M$.

In the Finsler space, the scalar product
 is essentially non-linear object with respect to
 the entered vectors, so that we may hope to meet  similar preservation property
in the Finslerian domains
if only  we  apply the connection which is non-linear, in the sense that
 the involved connection coefficients depend on tangent vectors
 $y$
 in non-linear way.
With this hope,  we need
 the {\it metrical  non-linear Finsler   connection}
 ${\cal FN}$,
 such that
 \be
{\cal FN}=\{N^m{}_j,D^m{}_{ij}\}: \qquad   N^m{}_j=N^m{}_j(x,y),
\qquad D^m{}_{ij}=D^m{}_{ij}(x,y).
\ee
The adjective ``metrical'' means that the action of the entailed
covariant derivative $\cD$ on the Finsler metric function,
 and also on the Finsler metric tensor, yields identically zero.
The coefficients $N^m{}_j$ and $D^m{}_{ij}$
are assumed to be positively homogeneous regarding the dependence on  vectors $y$,
respectively of degree 1 and degree 0.

Accordingly, the most important object what should be lifted from
the Riemannian to Finslerian
space
is
the two-vector angle,
to be denoted by
$ \al_{\{x\}}(y_1,y_2)$,
where
$
y_1,y_2\in T_xM.
$
Like to the Riemannian geometry proper,
the underlined idea is to measure the angle by means of length of
the respective geodesic arcs evidenced  on the indicatrix.

The Finsler
 space endows the vector pair $y_1,y_2$ with the scalar product
$$
{\langle y_1,y_2\rangle}_{\{x\}}=F(x,y_1)F(x,y_2)
\al_{\{x\}}(y_1,y_2)
$$
 on analogy of the Riemannian geometry.

The non-linear deformation
\be
{\cal FN}=\cC\cdot {\cal RL}
\ee
of the Riemannian connection
may exist to yield the Finsler connection
${\cal FN}$ which preserves the Finslerian
two-vector angle
$ \al_{\{x\}}(y_1,y_2)$
under the associated parallel transports of the  vectors
$y_1,y_2$.

In the theory of Finsler spaces, the key objects, the connection included,
 were introduced and studied
on the basis of various convenient sets of axioms
(see [1-5] and references therein).
Regarding the significance of the angle notion,
the important farther step was made in [6]
were in processes of studying  implications of the
two-vector angle defined by area,
 the theorem was proved
 which states that a diffeomorphism
between two Finsler spaces is an isometry iff it keeps the  angle.
 This  Tam{\' a}ssy's  theorem
substantiates
the  idea  to develop the Finsler connection from the
Finsler two-vector angle, possibly on the analogy of  the Riemannian geometry.

To meet new methods of applications,
the interesting chain of linear connections
was introduced and studied in [3].
It was emphasized that in the Riemannian geometry we have naturally
the metrical and linear connection.
We depart  from this connection
to develop the Finsler connection.

Namely,
we shall confine our attention to the   case when the  space ${\cal F}^N$
is obtainable from
the  space ${\cal R}^N$
by means of the {\it conformal automorphism},
according to the definition (2.1) of Section 2.
We shall also assume that
under the used transformations
the Finslerian indicatrix
$
{\cal  IF}_{\{x\}}\in T_{x}M
$
and
the Riemannian sphere
$
{\cal S}_{\{x\}}\in T_{x}M
$
are in correspondence (according to (2.2)).

Additionally, we subject the $\cC$-transformations
 to the condition of positive homogeneity
with respect to tangent vectors $y$, denoting the degree of homogeneity
 by $H$.
We call the $H$ the {\it degree of conformal automorphism}.

Remarkably,
such Finsler spaces of  dimensions $N\ge 3$
can be characterized by the condition that
the  indicatrix is a space of constant curvature
(see Proposition 2.1).
{\it The indicatrix curvature value
 ${\cal C}_{\text{Ind.}} $
is the square of
the degree of conformal automorphism,}
 that is,
\be
{\cal C}_{\text{Ind.}} \equiv H^2
\ee
 (indicated in  (2.3)).
 The condition has been realized, the Finslerian two-vector angle
 $ \al_{\{x\}}(y_1,y_2)$
proves to be a factor of
 the angle
$ \al^{\text{Riem}}_{\{x\}}(y_1,y_2) $
operative traditionally in the Riemannian space,
namely the simple equality
\be
  \al_{\{x\}}(y_1,y_2)=\fr1{H(x)}  \al^{\text{Riem}}_{\{x\}}(y_1,y_2)
\ee
(see (2.31)-(2.32))
is obtained.

The equality
\be
S(x,\bar y)=\lf(F(x,y)\rg)^H
\ee
(see (2.10)) is arisen,
which validates the indicatrix correspondence principle (2.2).

We  set forth the conventional requirement of preservation of the Finsler metric function $F(x,y)$, namely
\be
d_iF=0
\ee
with
\be
d_i=\D{}{x^i}
+N^k{}_{i}(x,y)
\D{}{y^k}.
\ee

With the definition
\be
\cD y^n~:=dy^n-N^n{}_j(x,y)dx^j
\ee
\ses\\
of covariant displacement of the  tangent vector,
the parallel transport of the   vector means
the  vanishing
\be
\cD y^n=0.
\ee

We apply this observation to the two-vector angle
$\al_{\{x\}}(y_1,y_2)$:
the   coefficients
$N^k{}_{i}(x,y)$
fulfill
the {\it angle preservation equation}
\ses\\
\be
d_i \al_{\{x\}}(y_1,y_2)=0,
\qquad y_1,y_2\in T_xM
\ee
under the parallel displacements of the entered vectors
$y_1$ and $y_2$,
if the involved operator $d_i$ is taken to read
\ses
\be
d_i=\D{}{x^i}
+N^k{}_{i}(x,y_1)
\D{}{y^k_1}
+
N^k{}_{i}(x,y_2)
\D{}{y^k_2}.
\ee
The $N^k{}_{i}(x,y)$ thus appeared can naturally be interpreted as the
{\it coefficients of the   non-linear connection produced by  angle}.

In this way we fulfill   the  canonical  geometrical principle:
the angle
$\al_{\{x\}}(y_1,y_2)$ formed by two vectors
$y_1$ and $y_2$
is left unchanged
under the parallel displacements of the  vectors
$y_1$ and $y_2$,
namely
$\cD\al\eqdef(dx^i)d_i\al=0$, for $d_i\al=0$.

In general the indicatrix curvature value
${\cal C}_{\text{Ind.}} $
may depend on the points $x\in M$.
We say that the space   ${\cal F}^N$
is {\it indicatrix-homogeneous}, if the value is a constant.
In view of the result ${\cal C}_{\text{Ind.}} \equiv H^2$
(indicated in  (2.3)),  such spaces can be characterized by the condition
that
the degree  $H$ of conformal automorphism
 is independent of $x$.

It proves that
{\it in the indicatrix-homogeneous case of the studied space}
 ${\cal F}^N$
 {\it  the equations}
(1.13)--(1.14)
 {\it can explicitly be solved for the coefficients}
 $N^k{}_{i}$
(see Proposition 2.2 and Note placed thereafter in Section 2).

From the obtained coefficients
$N^k{}_{m} $ given by (2.36),
the entailed coefficients
\be
N^k{}_{mn}=
\D{N^k{}_{m} }{y^n},  \qquad
N^k{}_{mnj}=
\D{N^k{}_{mn} }{y^j}
\ee
can straightforwardly be evaluated  (Section 3).
Let us use the coefficients to construct  the covariant derivative
$\cD_mg_{nj}$ of the Finsler metric tensor
$g_{nj}=g_{nj}(x,y)$ of the considered space
${\cal F}^N$, namely
\be
\cD_mg_{nj}~:= d_mg_{nj}+ N^k{}_{mj}g_{kn} + N^k{}_{mn}g_{kj},
\ee
where $d_m$ is given by (1.10).
It proves that the covariant derivative introduced by (1.16)
with  the  coefficients
$N^k{}_{m} $ given by (2.36)
 possesses the property
\be
\cD_mg_{nj}=0
\ee
in the indicatrix-homogeneous case.
 The property  can be verified by straightforward substitutions which result in
the vanishing
\be
y_kN^k{}_{mnj}=0
\ee
(see Proposition 3.1).

It is amazing but the fact that the last vanishing is an implication of the identity
$y^kC_{knj}=0$ shown by the   Cartan tensor  $C_{knj}=(1/2)\partial g_{kn}/\partial y^j$.
Indeed,
additional evaluation leads to
the result
\be
N^k{}_{mnj}=
-
\cD_mC^k{}_{nj}
\ee
in the
indicatrix-homogeneous case (see Proposition 3.2),
where
\be
\cD_mC^k{}_{nj}~:=
d_mC^k{}_{nj}-N^k{}_{mt}C^t{}_{nj}
+N^t{}_{mn}C^k{}_{tj}
+N^t{}_{mj}C^k{}_{nt}.
\ee

The coefficients
$
N_{kmnj}
=
g_{kh}N^h{}_{mnj}
$
can be written as
\be
N_{kmnj}=
-
\cD_mC_{knj}
\ee
and, therefore, they are {\it symmetric}
with respect to the subscripts $k,n,j$.

Thus, with the identification
\be
D^k{}_{in}(x,y)=-N^k{}_{in}(x,y),
\ee
 in    the  Finsler space
${\cal F}^N$  of the indicatrix-homogeneous type (that is, when $H=const$)
the metrical  angle-preserving connection (1.4)
is given by the coefficients
 $\{N^k_i, D^k{}_{in}\}$
 found explicitly.
 Recollecting the assumed homogeneity of the coefficients,
 from (1.22) we infer the equality
\be
D^k{}_{im}y^m=-N^k{}_{i}.
\ee

Realizing the $\cC$-transformation
locally by
$y^i=y^i(x,t)$ with $ t^n \equiv \bar y^n
$
(see (2.11))
and applying
 the Riemannian operator
$$
d^{\text{Riem}}_i=
\D{}{x^i}-a^k{}_{ih}t^h\D{}{t^k}
$$
(cf. (1.3)) to the field $y^i(x,t)$,
it is possible to conclude that
\be
N^n{}_{i}=d^{\text{Riem}}_i   y^n
\ee
(see (2.47)).
This representation of the coefficients $N^n{}_{i}$
possesses a clear geometrical and tensorial meaning
and is  alternative (and equivalent)
to the representation
(2.36).
The derivation of the representation (1.24) uses the formula (1.23).

According to Proposition 2.3,
   the  Finsler space
${\cal F}^N$  of the indicatrix-homogeneous type
is obtained from
the Riemannian space
${\cal R}^N$
by means of the {\it parallel deformation}.

Since also the transitivity of covariant derivative holds, namely
$
\cD_nt^i=0
$
(see (2.39)), and
$ g_{kh}= {\cal C}^m_k   {\cal C}^n_h a_{mn}$
(see (2.25)),
we should conclude that
in    the  Finsler space
${\cal F}^N$  of the indicatrix-homogeneous type
 the metrical angle-preserving connection
 is  the $\cC$-export of
 the
  metrical  linear Riemannian  connection
(1.2) applied conventionally
in
the background  Riemannian space
 ${\cal R}^N$.

In  Section  4 we perform the attentive
comparison between the commutators
of the involved Finsler  covariant derivative $\cD$ and
 the commutators of  the underlined
Riemannian
covariant derivative $\nabla$,
 assuming $H=const$.
By this method,  we derive
the associated curvature tensor
 $\Rho_k{}^n{}_{ij}$.
 The skew-symmetry
$  \Rho_{mnij}=-\Rho_{nmij}=- \Rho_{mnji}
$
 holds.
   The covariant derivative
 $\cD_l$ of the tensor fulfills
  the cyclic identity,
  completely similar to the  Riemannian case
in which
  the cyclic identity is valid for the derivative  $\nabla_la_k{}^n{}_{ij}$
of the
 Riemannian curvature tensor
$a_k{}^n{}_{ij}$.
The tensor
$
 M^n{}_{ij}=- y^k\Rho_k{}^n{}_{ij}
$
proves to be transitive to the Riemannian tensor
$-y^n_t t^h
a_h{}^t{}_{ij},
$
namely the equality
$
M^n{}_{ij}=
-
y^n_t t^h
a_h{}^t{}_{ij}
$
holds.
The very tensor $\Rho_k{}^n{}_{ij}$ is not transitive
to the Riemannian precursor
$a_h{}^m{}_{ij} $, instead the more general equality
$$
\Rho_k{}^n{}_{ij}
=
-
(1-H)\fr1{F}(l_k\de^n_m-l^ng_{mk})
M^m{}_{ij}
 +
y^n_ma_h{}^m{}_{ij} t^h_k
$$
is obtained.
We observe that the difference between the curvature tensor
$\Rho_k{}^n{}_{ij}$
and the transitive term
$y^n_ma_h{}^m{}_{ij} t^h_k$
is proportional to
$(1-H)$.
Squaring the tensor yields the sum
$$
\Rho^{knij}\Rho_{knij}=
a^{knij}a_{knij}
+\fr 2{S^2}
\lf(\fr1{H^2}-1\rg)
t^l
a_l{}^{nij}
t^h
a_{hnij},
$$
\ses\\
which is the  ${\cal F}^N$-extension of the Riemannian term
$a^{knij}a_{knij}.$ The difference
$ \Rho^{knij}\Rho_{knij}  -  a^{knij}a_{knij}$  is proportional to $(H^{-2}-1)$.

In Section  5 we develop an explicit and attractive particular case, namely we present
the explicit example (5.27) of the conformally automorphic transformation (2.1),
specializing the Finsler space to be the ${\cal FS}$-space.
The  space is endowed with the Finsler metric function $F$
which is constructed from a Riemannian metric function
 $
 S= \sqrt{a_{ij}(x)y^iy^j}
$
and an 1-form
$ b= b_i(x)y^i$
according to
 the functional dependence
\be
F(x,y) =\Phi \lf(x; b, S,y\rg),
\ee
where $\Phi $ is a sufficiently smooth scalar function.
In step-by-step way, we derive  the  coefficients $N^m{}_i$ specified by (2.36),
obtaining  the  explicit representation (5.72)--(5.75).
It proves that the suitability of the proposed transformation
(5.27) imposes the severe restriction on the Finsler metric function,
namely the function  must be of the Finsleroid type
(described in [7]).
In the restricted case which implies independence
of   the function
$\Phi \lf(x; b, S,y\rg)$ of $x$, assuming also that the Riemannian norm
of the 1-form $b$ is a constant,
the obtained   coefficients $N^m{}_i$
straightforwardly entail the  vanishing set
$
\cD_n F=\cD_n y_j=\cD_n g_{ij}=0
$
(see (5.96)-(5.98)),
together with the angle preservation (1.13).
Simplifying coefficients $N^m{}_i$
culminates in the representation (5.102).
The initial transformation
(5.27)
reduces to (5.106),
for it proves possible to find explicitly the involved functions
$\varrho$ and $\mu$.

In Conclusions, Section 6, we emphasize several important ideas.

In Appendices A-E
we present the explicit evaluations which are required to verify the validity
of the formulated propositions.

\ses

\ses

\setcounter{sctn}{2} \setcounter{equation}{0}

\bc  {\bf 2. Main observations}  \ec

\ses

\ses

Below, {\it  any dimension }  $N\ge 3$ is allowable.

\ses

Let $M$ be an $N$-dimensional
$C^{\infty}$
differentiable  manifold, $ T_xM$ denote the tangent space to $M$ at a point $x\in M$,
and $y\in T_xM\backslash 0$  mean tangent vectors.
Suppose we are given on the tangent bundle $TM$ a Riemannian metric
 $S$.
 Denote by
$\cR^N=(M,S)$
the obtained $N$-dimensional Riemannian space.
Let additionally a Finsler metric function  $F$ be introduced on this  $TM$,
yielding
a Finsler space
 ${\cal F}^N=(M,F)$.
We shall study the Finsler space
 ${\cal F}^N$
 specified according to the following definition.

\ses

\ses

INPUT DEFINITION.
The  space ${\cal F}^N$ is  {\it conformally   automorphic} to
the Riemannian space
${\cal R}^N$:
\be
{\cal F}^N=\cC\cdot {\cal R}^N
\ee
such that in each tangent space $T_xM$ the $\cC$-automorphism
transforms  conformally the metric produced by the Finsler metric  to
the Euclidean metric entailed by the Riemannian metric.
It is assumed that  the applied $\cC$-transformations
 do not influence any point $x\in M$ of the base manifold  $M$
 and that they are invertible.
It is also natural to require that
the $\cC$-transformations send  unit vectors to unit vectors:
\be
{\cal  IF}_{\{x\}}= \cC \cdot    {\cal S}_{\{x\}}.
\ee
Additionally, we subject the $\cC$-transformations
 to the condition of positive homogeneity
with respect to tangent vectors $y$, denoting the degree of homogeneity
 by $H$.
 We call the $H$ the {\it degree of conformal automorphism}.

The existence of such spaces is explained by the following proposition.

\ses

\ses

{\bf Proposition 2.1.}
{\it
A Finsler space  is of the claimed type
${\cal F}^N$ if and only if the indicatrix of the Finsler space
 is a space of constant curvature. Denoting the indicatrix curvature value by
${\cal C}_{\text{Ind.}}$,
the equality
\be
{\cal C}_{\text{Ind.}} \equiv H^2, \qquad H>0,
\ee
is obtained.
The relevant  conformal multiplier is given by $p^2$ with
\be
p=
\fr1{H}
F^{1-H}.
\ee

}

\ses

\ses

The proposition is of the local meaning in both the base manifold and the
tangent space.
The validity of the proposition can be verified by simple straightforward evaluations
(which are presented in Appendix A).

The value
${\cal C}_{\text{Ind.}}$
may vary from point to point of the manifold
$M$, so that in general $H=H(x)$.

We take ${\cal C}_{\text{Ind.}} >0$.
 Extension of the proposition to negative value
of
${\cal C}_{\text{Ind.}}$ would be a straightforward task.

On every punctured tangent space $T_xM\setminus0$,
the Finsler metric function $F$ is assumed to be positive,
 and also  positively homogeneous
of degree 1:
\be
F(x,ky)=kF(x,y), \qquad k>0, \forall y.
\ee
 Therefore,  the conformal factor $p^2=\lf(F^{1-H}/H\rg)^2$
possesses this kind of homogeneity with degree $2(1-H)$.
For a given  function $F$
we  can construct
the covariant tangent vector $\hat y=\{y_i\}$
and the  Finslerian metric tensor $\{g_{ij}\}$
in the ordinary way:
$$
y_i :=\fr12\D{F^2}{y^i}, \qquad
g_{ij} :
=
\fr12\,
\fr{\prtl^2F^2}{\prtl y^i\prtl y^j}
=\fr{\prtl y_i}{\prtl y^j}.
$$
\ses
The contravariant tensor $\{g^{ij}\}$ defined by the reciprocity conditions
$g_{ij}g^{jk}=\de^k_i$,
where $\de$ stands for the Kronecker symbol.

Let the $\cC$-transformation (2.1) be assigned locally by means of the differentiable functions
\be
\bar y^m=\bar y^m(x,y),
\ee
subject to the required   homogeneity
\be
\bar y^m(x,ky)=k^H\bar y^m(x,y), \qquad k>0, \forall y.
\ee
This entails the identity
\be
\bar y^m_ky^k=H\bar y^m,
\ee
where
$  \bar y^m_k=\partial \bar y^m/\partial y^k$.
Fulfilling  the conformal automorphism (2.1) means locally
\be
g_{mn}(x,y)=c_{ij}(x,\bar y) \bar y^i_m   \bar y^j_n, \qquad  c_{ij}=p^2a_{ij}(x).
\ee

\ses

Contracting the $g_{mn}$ by $y^my^n$ and noting
the
involved  homogeneity
together with the value (2.4) of $p$, we get the equality
\be
S(x,\bar y)=\lf(F(x,y)\rg)^H,
\ee
where
$ S=\sqrt{a_{mn}(x)\bar y^m\bar y^n}. $

\ses

Denote
by
\be
y^i=y^i(x,t), \qquad t^n \equiv \bar y^n,
\ee
the inverse transformation, so that
\be
y^i(x,k t)=k^{1/H} y^i(x,t), \qquad k>0, \forall t,
\ee
and
\be
y^i_nt^n=\fr1Hy^i,
\ee
where
$ y^i_n=\partial y^i/\partial t^n$.
The inverse to (2.9) reads:
\be
g_{kh} y^k_my^h_n  =  c_{mn}.
\ee

\ses

The following useful relations can readily be arrived at:
\ses\\
\be
y_my^m_n=\fr{F^2}{HS^2}t_n
\equiv
\fr1HF^{2(1-H)}t_n, \qquad t_n=a_{nh}t^h,
\ee
and
\ses\\
$$
y_my^m_{nl}t^l_j
+
g_{mj}y^m_n
=
2
\lf(\fr1H-1\rg)
F^{-2H}y_jt_n
+
\fr1HF^{2(1-H)}a_{nh}t^h_j,
$$
\ses
where
$t^l_j=\bar y^l_j,
y^m_{nl}=\partial { y^m_n}/\partial{ y^l}.$
Alternatively,
\be
t_ht^h_n=\fr{HS^2}{F^2}y_n
\equiv
HF^{2(H-1)}y_n
\ee
\ses\\
and
\be
t_ht^h_{nu}y^u_i
+
a_{hi}t^h_n
=
2
(H-1)
F^{-2}t_iy_n
+
HF^{2(H-1)}g_{nu}y^u_i,
\ee
where $t^h_{nu}=\partial t^h_{n}/\partial y^u.$
We may  also write
\ses\\
$$
t_ht^h_{ni}
+
H^2  F^{2(H-1)}g_{ni}
=
2
(H-1)
F^{-2}HF^{2(H-1)}y_iy_n
+
HF^{2(H-1)}g_{ni},
$$
or
\ses\\
\be
t_ht^h_{ni}
=
H(1-H)F^{2(H-1)}(g_{ni}-2l_nl_i).
\ee

\ses

From (2.14) it follows that
$
g_{nm}y^m_i=p^2t^j_n a_{ij}.
$

Differentiating (2.9) with respect to $y^k$ yields the following
representation for the Cartan tensor
$C_{mnk}=(1/2)\partial g_{mn}/\partial y^k$:
\be
2C_{mnk}=
(1-H)\fr2{F}l_k
g_{mn}
+p^2 (t^i_{mk}t^j_n +t^i_mt^j_{nk}) a_{ij}.
\ee
Contracting this tensor by $y^n$ results in the equality

\be
p^2 t^i_{mk}t^j  a_{ij}=
\lf(\fr1H-1\rg)(h_{km} -l_k  l_m ),
\ee
where the vanishing $C_{mnk}y^n=0$ and the homogeneity identity (2.8)
have
been taken into account.

\ses

Symmetry of the tensor $C_{mnk}$ demands

\be
(1-H)\fr2{F}(l_kg_{mn}-l_mg_{kn})
+p^2 (t^i_mt^j_{nk}-t^i_kt^j_{nm}) a_{ij}
=0,
\ee
so that we may alternatively write
\be
C_{mnk}=
(1-H)\fr1{F}(l_kg_{mn}+l_ng_{mk}-l_mg_{nk})
+p^2 t^i_mt^j_{nk} a_{ij}.
\ee

Contracting  the last tensor
by $g^{nk}$ yields
\ses\\
$$
2C_{m}=
(1-H)\fr2{F}l_m
+g^{nk}
p^2 (t^i_{nk}t^j_m +t^i_nt^j_{mk}) a_{ij}
\equiv 2C_{mnk}g^{nk},
$$
from which it ensues that
$$
2C_{m}=
(1-H)\fr2{F}l_m
+
2
g^{nk}
p^2 t^i_{nk}t^j_m  a_{ij}
+
g^{nk}
p^2 (t^i_nt^j_{mk}-t^i_mt^j_{nk})
 a_{ij},
$$
or
\ses\\
$$
2C_{m}=
(1-H)\fr2{F}l_m
+
2
g^{nk}
p^2 t^i_{nk}t^j_m  a_{ij}
-
(1-H)
g^{nk}
\fr2{F}(l_mg_{nk}-l_ng_{mk})
 a_{ij}.
$$
It is also convenient to use the representation
\ses\\
\be
FC_{m}=
-(N-2)(1-H)l_m
+
Fg^{nk}
p^2 t^i_{nk}t^j_m  a_{ij}.
 \ee

 The  space ${\cal F}^N$ is obtainable
  from
  the Riemannian
 space
 ${\cal R}^N$
  by means of the  deformation (1.1)
which, owing to  (2.2) and (2.9),
can be presented by the {\it conformal deformation tensor}
\be
{\cal C}^m_k~:=p\bar y^m_k,
\ee
so that
\be
 g_{kh}= {\cal C}^m_k   {\cal C}^n_h a_{mn}.
\ee
\ses
The zero-degree homogeneity
\be
{\cal C}^m_n(x,ky)=  {\cal C}^m_n(x,y), \qquad k>0,\forall y,
\ee
holds,
together with
\be
{\cal C}^m_ny^n=F^{1-H}\bar y^m.
\ee
\ses

The indicatrix  correspondence  (2.2) is a direct implication of the
equality $S=F^H$ (see (2.10)).
We may apply the considered transformation (2.6) to the unit vectors:
\be
l=\cC\cdot L: ~ ~ ~  l^i=y^i(x,L);
~ \qquad
L=\cC^{-1}\cdot l: ~ ~  ~  L^i=t^i (x,l),
\ee
where
$l^i=y^i/F(x,y)$ and $L^i=t^i/S(x,t)$ are components of the respective
Finslerian and Riemannian unit vectors,
which possess the properties $F(x,l)=1$ and $S(x,L)=1$.
We have $L^m=t^m(x,l)$.
On the other hand, from (2.4) and (2.9) it just follows that
\be
g_{mn}(x,l)=\fr1{H^2}a_{ij}(x)t^i_m(x,l) t^j_n(x,l),
\ee
\ses
so that  under the transformation (2.28) we have
\be
g_{mn}(x,l)dl^mdl^n=\fr1{H^2}a_{ij}(x)dL^i dL^j.
\ee

No support vector enters  the right-hand part in the previous  equality
(2.30).
Therefore,
any two nonzero tangent vectors
  $y_1,y_2\in T_xM$ in  a fixed tangent space  $T_xM$
form the  ${\cal F}^N${\it-space angle}
  \be
 \al_{\{x\}}(y_1,y_2) =
 \fr1{H(x)}
 \arccos\la,
  \ee
 where the scalar
\be
\la =
  \fr{a_{mn}(x)t_1^m  t_2^n}
 { S_1S_2 },
\qquad \text{with} \quad  t_1^m = t^m (x,y_1)  \quad   \text{and} \quad
t_2^m = t^m (x,y_2)    ,
   \ee
is of the entire Riemannian  meaning in the
space ${\cal R}^N$;
the notation
$$
S_1=\sqrt{a_{mn}(x)t^m_1t^n_1}, \qquad S_2=\sqrt{a_{mn}(x)t^m_2t^n_2}
$$
has been used.

From (2.32) it follows that
$$
\D{\la}{x^i}=
  \fr{a_{mn,i}t_{1}^m  t_2^n}
{S_1S_2 }+
\fr1{ S_1S_2}
  a_{mn}\lf(\D{t_{1}^m}{x^i}  t_2^n
  +
t_{1}^m  \D{t_2^n  }{x^i}
\rg)
$$

\ses

\be
-
\fr12\la
\lf[
\fr1{ S_1S_1 }
\lf(a_{mn,i}t_{1}^mt_1^n
+ 2 a_{mn} \D{t_{1}^m}{x^i} t_1^n
 \rg)
    +
\fr1{S_2S_2 }
\lf(a_{mn,i}t_{2}^mt_2^n
+ 2 a_{mn} \D{t_{2}^m}{x^i}  t_2^n
 \rg)
\rg],
\ee
\ses
where   $a_{mn,i}=\partial a_{mn}/\partial x^i$,
and
\be
\D{\la}{y^k_1}=
\lf[
  \fr{a_{mn}  t_2^n}
{S_1S_2 }
 -      \fr{a_{mn} t_1^n}
{S_1S_1}\la
\rg]t_{1k}^m ,
\qquad
\D{\la}{y^k_2}=    \lf[
 \fr{a_{mn}  t_1^n}
{ S_2S_1 }
 -      \fr{a_{mn} t_2^n}
{S_2S_2}\la
\rg]t_{2k}^m.
\ee

Let  the coefficients $N^k{}_i$ be subjected to the equation
\ses\\
\be
d_i\la=0,
\ee
where $d_i$ is the operator (1.10).

It is possible to establish the validity of the following proposition.

\ses

\ses

{\bf Proposition 2.2.}
{\it
 When} $d_iF=0$  {\it and} $H=const$,
 {\it the equation} (2.35)
  {\it can  be solved for the coefficients
$N^m{}_n$,
yielding
\ses\\
\be
N^m{}_n=-y^m_i \lf(\D{t^i}{x^n}+a^i{}_{kn}t^k\rg).
\ee

}

\ses

\ses

See Appendix B.

\ses

In (2.36), the $a^i{}_{kn}=a^i{}_{kn}(x)$ are
the Christoffel symbols
\be
a^i{}_{kn}
=
\fr12
a^{ih}\Bigl( \D{a_{hk}}{x^n}+ \D{a_{hn}}{x^k}- \D{a_{kn}}{x^h} \Bigr)
\ee
of the Riemannian space
$\cR^N$.

\ses

\ses

{\bf Note.}   When $H=const$,  from (2.31)   it just follows that the angle
$ \al_{\{x\}}(y_1,y_2)$
fulfills the vanishing which is completely similar to
(2.35),   namely the vanishing  (1.13)  claimed in Section 1.

\ses

\ses

With the covariant derivative
\be
\cD_nt^i~:=d_nt^i+a^i{}_{kn}t^k
\ee
the representation (2.36) can be interpreted as the
manifestation of the
{\it transitivity}
\be
\cD_nt^i=0
\ee
of the  connection
 under the  conformal   automorphism (2.1).

By differentiating (2.39) with respect to $y^m$
 we may conclude that
 the covariant derivative
\be
\cD_n t^i_m ~:= d_n t^i_m-D^h{}_{nm}t_h^i+a^i{}_{nl}t^l_m, \qquad D^h{}_{nm}=-N^h{}_{nm},
\ee
vanishes identically:
\be
\cD_n t^i_m= 0.
\ee

Since
$ y^n_k t^k_j=\de^n_j$,
the previous identity can be transformed to
\be
d^{\text{Riem}}_i y^n_k+D^n{}_{is}y^s_k-a^h{}_{ik}y^n_h=0,
\ee
which can be interpreted as the  covatiant derivative
vanishing:
\be
\cD_i y^n_k=0.
\ee
This formula entails
\be
\cD_i y^n=0
\ee
(because of (2.39)), where
\be
\cD_i y^n~:=
d^{\text{Riem}}_i
 y^n
 +D^n{}_{is}y^s.
\ee
Here, $y^n$ mean the functions $y^n(x,t)$   introduced by (2.11).
We  have
used the Riemannian operator
\be
d^{\text{Riem}}_i=
\D{}{x^i}+L^k{}_i\D{}{t^k}, \qquad
 L^k{}_i=   -a^k{}_{ih}t^h
\ee
(cf. (1.3)).

Since $D^n{}_{is}y^s=-N^n{}_{i}$, from (2.44)-(2.45) we may conclude that the representation
\be
N^n{}_{i}=d^{\text{Riem}}_i   y^n
\equiv
 \D{y^n(x,t)}{x^i}  +   y^n_hL^h{}_i
  \ee
is valid which is alternative to (2.36).

Let us verify (2.42). We have
$$
0= y^n_k
 \lf(
\D{t^k_j}{x^i}    + N^h{}_it^k_{hj}-D^h{}_{ij}t_h^k+a^k{}_{il}t^l_j
\rg)
 $$

\ses

\ses

$$
=
-t^k_j \lf(
\D{ y^n_k}{x^i}  + y^n_{kh} \D{ t^h}{x^i}\rg)
+
y^n_k
\lf(
N^h{}_it^k_{hj}-D^h{}_{ij}t_h^k+a^k{}_{il}t^l_j
 \rg).
 $$
Contracting this  by $y^j_m$ yields
\ses\\
$$
0=
\D{ y^n_m}{x^i}  + y^n_{mh} \D{ t^h}{x^i}
-
y^j_my^n_k
N^h{}_it^k_{hj}
+D^n{}_{ij}y^j_m-y^n_ka^k{}_{im}.
 $$
Take
$ N^h{}_i$
from (2.36):
$$
0=
\D{ y^n_m}{x^i}  + y^n_{mh} \D{ t^h}{x^i}
+
y^j_my^n_k
y^h_l
t^k_{hj}
 \lf(\D{t^l}{x^i}+a^l{}_{ji}t^j\rg)
+D^n{}_{ij}y^j_m-y^n_ka^k{}_{im}.
 $$
This vanishing is tantamount to the considered (2.42), for
$
y^j_my^n_k
y^h_l
t^k_{hj}
=
- y^n_{ml}.
$

Owing to  the equalities (2.4), (2.24), and (2.43),
we are entitled to formulate the following proposition.

 \ses

\ses

{\bf Proposition 2.3.}
{\it
 When} $d_iF=0$  {\it and} $H=const$,
 {\it
the  deformation tensor} (2.24) {\it is  parallel}
\be
\cD_n {\cal C}^m_k=0,
\ee
{\it where}
\be
\cD_n {\cal C}^m_k= d_n {\cal C}^m_k -D^h{}_{nk}{\cal C}_h^m+a^m{}_{nl}{\cal C}^l_k.
\ee

\ses

\ses

With these observations,
it is possible
 to develop a direct method to induce the
connection   in the Finsler
 space    ${\cal F}^N$ from the
  metrical  linear Riemannian  connection  (1.2)
meaningful in the background Riemannian space ${\cal R}^N$.

The coefficients $N^k{}_i(x,y)$
can also be obtained by means of
 the
transitivity map
\be
\{N^k{}_i\}=\cC\cdot \{L^k{}_i\}.
\ee

Indeed, with
an arbitrary  differentiable scalar
$w(x,y)$,
we can apply the transformation
$\{y^i=y^i(x,t),\, t^n \equiv \bar y^n\}$
indicated in (2.11) and consider the $\cC$-transform
\be
W(x,t)=w(x,y), \quad \text{ which entails} \quad \D W{t^n}= y^k_n \D w{y^k},
\ee
thereafter   postulating  that the $\cC$-transformation is {\it covariantly transitive},
namely
\be
\Biggl(\fr{\partial}{\partial x^i}
+N^k{}_i{(x,y)} \fr{\partial}  {\partial y^k}
\Biggr)
w(x,y)
=
\Biggl(\fr{\partial}{\partial x^i}
+L^k{}_i{(x,t)} \fr{\partial}{\partial t^k}
\Biggr)
W(x,t).
\ee
Since the field $w$ is arbitrary, the last equality is fulfilled if and only if
\be
 N^k{}_i=d^{\text{Riem}}_iy^k
\equiv
 \D{y^k(x,t)}{x^i}  +   y^k_hL^h{}_i.
\ee
This is the    representation  which is required  to
realize the map (2.50).
We have again arrived at the coefficients (2.47).

With the knowledge of the coefficients
$N^k{}_i(x,y)$,
 we can use the formulas (2.40) and (2.41)
to express the Finslerian connection coefficients
$D^h{}_{nm}$
through
the Riemannian Christoffel symbols
$a^i{}_{nl}$.
Thus
we have   induced the
connection   in the Finsler
 space    ${ \mathcal   F}^N$ from the
  metrical  linear Riemannian  connection  (1.2)
meaningful in the background Riemannian space ${\mathcal   R}^N$.

It can readily be noted that the transitivity property (2.52) can
straightforwardly
be extended to scalars
dependent on two vectors. Namely,
if
\be
W(x,t_1,t_2)=w(x,y_1,y_2),
\ee
then
\be
\Biggl(\fr{\partial}{\partial x^i}
+N_{1i}^k \fr{\partial}  {\partial y_1^k}
+N_{2i}^k \fr{\partial}  {\partial y_2^k}
\Biggr)
w(x,y_1,y_2)
=
\Biggl(\fr{\partial}{\partial x^i}
+ L_{1i}^k \fr{\partial}{\partial t_1^k}
+ L_{2i}^k \fr{\partial}{\partial t_2^k}
\Biggr)
W(x,t_1,t_2),
\ee
where $N_{1i}^k=N^k{}_i(x,y_1),\,N_{2i}^k=N^k{}_i(x,y_2),\,
L_{1i}^k =L^k{}_i{(x,t_1)}, \,L_{2i}^k =L^k{}_i{(x,t_2)} $.

\ses

\ses

\setcounter{sctn}{3} \setcounter{equation}{0}

\bc  {\bf 3. Properties of  connection coefficients}\ec

\ses

\ses

The derivative coefficients (1.15)
can straightforwardly be evaluated from (2.36).
We obtain explicitly
\be
N^k{}_{mn}
=
   -  y^k_{sl}t^l_n T^s_m     -  y^k_s T^s_{n,m}   \quad \text{with}
\qquad
T^s_m =        \D{t^s}{x^m}  +a^s{}_{mh}t^h , \quad
 T^s_{n,m}= \D{t^s_n}{x^m}    +   a^s{}_{mh}t^h_n   ,
\ee
\ses\\
which entails the contractions
\be
y_kN^k{}_{mn}
=
-
\lf(
2
\lf(\fr1H-1\rg)
F^{-2H}y_nt_s
+
\fr1HF^{2(1-H)}a_{sh}t^h_n
-
g_{ln}y^l_s
\rg)
T^s_m
   -
\fr{F^2}{HS^2}t_s
 T^s_{n,m}
 \ee
and
$$
y_kN^k{}_{mn}
+
g_{ln}N^l{}_m
=
-
\lf(
2
\lf(\fr1H-1\rg)
F^{-2H}y_nt_s
+
\fr1HF^{2(1-H)}a_{sh}t^h_n
\rg)
T^s_m
   -
\fr{F^2}{HS^2}t_s
 T^s_{n,m},
  $$
\ses\\
together with
$$
y_kN^k{}_{mni} +  g_{ki}N^k{}_{mn}
+
g_{ln}N^l{}_{mi}
+
2C_{lni}N^l{}_m
$$

\ses

\ses

$$
=
-
\lf(\fr1H-1\rg)
2F^{-2H}
\lf[
    (g_{ni}-2Hl_nl_i)   t_s
+
(y_na_{sl}t^l_i+y_ia_{sl}t^l_n)
\rg]
T^s_m
-
\fr1H
F^{2(1-H)}a_{sh}t^h_{ni}
T^s_m
$$

\ses

\ses

\ses

$$
-
\lf(
2
\lf(\fr1H-1\rg)
F^{-2H}y_nt_s
+
\fr1HF^{2(1-H)}a_{sh}t^h_n
\rg)
 T^s_{i,m}
 $$

          \ses

\ses

\be
   -
\lf(
2
\lf(\fr1H-1\rg)
F^{-2H}y_it_s
+
\fr1HF^{2(1-H)}a_{sh}t^h_i
\rg)
 T^s_{n,m}
   -
\fr1H
F^{2(1-H)}t_s
\lf(    \D{t^s_{ni}}{x^m}    +   a^s{}_{mh}t^h_{ni}   \rg).
\ee
\ses\\
The attentive calculation of the entered terms
(carried out in Appendix C)
leads to the following remarkable result.

\ses

\ses

{\bf Proposition 3.1.}
{\it
If the coefficients $N^k{}_m$ are taking according to
Proposition}
2.2,
{ then
the vanishing
}
$
y_kN^k{}_{mnj}=0
$
holds identically.

\ses

\ses

In performing involved calculation it is necessary to note
that
in view of  (2.15) and (2.36), we can write
$$
d_mF=
\D F{x^m}+N^k{}_n\D F{y^k} =  \D F{x^m}+N^k{}_nl_k
=\D F{x^m}-
\fr 1{FH}F^{2(1-H)}
t_s
T^s_m
$$
so that, because of $d_mF=0$,
the equality
\be
\D F{x^m}
=
\fr 1{FH}F^{2(1-H)}
t_s
T^s_m
\ee
is valid.

It is also possible to evaluate the covariant derivative
$\cD_mC^k{}_{nj}$
(see (1.20)),
using the equality
$
d_mg_{hn}=-N^t{}_{mh}g_{tn}  -N^t{}_{mn}g_{th}
$
entailed by the metricity (1.16).
This way leads to the following result.

\ses

\ses

{\bf Proposition 3.2.}
{\it
The representation
$
N^k{}_{mnj}=
-
\cD_mC^k{}_{nj}
$
is valid,
whenever
$d_mF=0$
and
$H=const$.
}

\ses

\ses

Proof of this proposition can be arrived at during a long chain of straightforward
substitutions (see Appendix D).

\ses

\ses

\setcounter{sctn}{4}
\setcounter{equation}{0}

\ses

\ses

\bc   {\bf  4. Entailed curvature tensor} \ec

\ses

\ses

Throughout the present section we assume that  $H=const$.
Given a tensor $w^n{}_k=w^n{}_k(x,y)$ of the tensorial type (1,1),
commuting the covariant derivative
\be
\cD_iw^n{}_k~:=   d_iw^n{}_k+D^n{}_{ih}w^h{}_k-D^h{}_{ik}w^n{}_h
\ee
 yields the equality
\be
\lf[\cD_i\cD_j-\cD_j\cD_i\rg] w^n{}_k=M^h{}_{ij}\D{w^n{}_k}{y^h}  -E_k{}^h{}_{ij}w^n{}_h
+E_h{}^n{}_{ij}w^h{}_k
\ee
with the  tensors
\be
 M^n{}_{ij}~:= d_iN^n{}_j-d_jN^n{}_i
\ee
and
\be
E_k{}^n{}_{ij}~: =
d_iD^n{}_{jk}-d_jD^n{}_{ik}+D^m{}_{jk}D^n{}_{im}-D^m{}_{ik}D^n{}_{jm}.
\ee
\ses

When the choice
$
 D^k{}_{in} =-N^k{}_{in}
$
is made (cf. (1.23)),
 the tensor (4.3) can be written in the form
\be
 M^n{}_{ij}= \D{N^n{}_j}{x^i}-\D{N^n{}_i}{x^j}
-N^h{}_iD^n{}_{jh}+N^h{}_jD^n{}_{ih}.
\ee
By applying the commutation rule (4.2)
to the particular choices $\{F, y^n, y_k,g_{nk}\}$
and noting the vanishing
$\{ {\mathcal D}_i F =  {\mathcal D}_i y^n=  {\mathcal D}_i y_k = {\mathcal D}_ig_{nk}=0\}$,
we obtain the identities
\be
y_n M^n{}_{ij}=0,
\qquad
y^kE_k{}^n{}_{ij}  = -   M^n{}_{ij},   \qquad
y_nE_k{}^n{}_{ij}  = M_{kij},
\ee
and
\be
E_{mnij}+E_{nmij}=2C_{mnh}
M^h{}_{ij} \quad \text{with} \quad C_{mnh}=\fr12\D{g_{mn}}{y^h}.
\ee
Differentiating (4.5) with respect to $y^k$ and using the equality
$
N^j{}_i=-D^j{}_{ik}y^k
$
 yield
\be
E_k{}^n{}_{ij}= -\D{M^n{}_{ij}}{y^k}.
\ee

The  cyclic identity
\be
\cD_kM^n{}_{ij}+\cD_jM^n{}_{ki}+\cD_iM^n{}_{jk}=0
\ee
is valid,
where
\be
\cD_k M^n{}_{ij}=d_k M^n{}_{ij}
+D^n{}_{kt} M^t{}_{ij}
 - a^s{}_{ki}   M^n{}_{sj}  -  a^s{}_{kj}   M^n{}_{is}.
\ee

It proves pertinent to replace in the commutator (4.2)
the partial derivative
$\partial  w^n{}_k/\partial y^h$
by the definition
\be
\cS_h{w^n{}_k}~:=\D{w^n{}_k}{y^h}+C^n{}_{hk}w^h{}_k-C^m{}_{hk}w^n{}_m
\ee
which has the meaning of the covariant derivative in the tangent
  space
supported by the point $x\in M$.
In particular,
$$
\mathcal S_h g_{nk}~:=\D{g_{nk}}{y^h} -C^m{}_{hn}g_{mk} -C^m{}_{hk}g_{nm}=0.
$$
With
the {\it curvature tensor}
\be
\Rho_k{}^n{}_{ij}~:=E_k{}^n{}_{ij}
-
M^h{}_{ij}C^n{}_{hk},
\ee
the commutator (4.2) takes on the form
\be
\lf(\cD_i\cD_j-\cD_j\cD_i\rg) w^n{}_k=
M^h{}_{ij}\cS_h w^n{}_k  -\Rho_k{}^h{}_{ij}w^n{}_h
+\Rho_h{}^n{}_{ij}w^h{}_k.
\ee

We denote
$\Rho_{knij}=
g_{mn}\Rho_k{}^m{}_{ij}.$
The skew-symmetry
\be
 \Rho_{mnij}=-\Rho_{nmij}
\ee
 holds (cf.  (4.7)).
The tensor obeys  also
the   cyclic identity
\be
\cD_l\Rho_k{}^n{}_{ij}+\cD_j\Rho_k{}^n{}_{li} +\cD_i\Rho_k{}^n{}_{jl} =0,
\ee
\ses
where
$$
\cD_l \Rho_k{}^n{}_{ij}
=
d_l\Rho_k{}^n{}_{ij}
+
D^n{}_{lt} \Rho_k{}^t{}_{ij}
-
D^t{}_{lk} \Rho_t{}^n{}_{ij}
 - a^s{}_{li}
 \Rho_k{}^n{}_{sj}
 -  a^s{}_{lj}
 \Rho_k{}^n{}_{is}.
$$

\ses

Let us realize the action of the $\cC$-transformation (2.1)-(2.2)
on tensors  by the help of the {\it transitivity  rule,}
that is,
\be
\{w^n{}_m(x,y)\}  =\cC\cdot \{W^n{}_m(x,t)\}: ~ ~ ~  w^n{}_m=  y^n_h t^j_mW^h{}_j,
\ee
where
 $W^n{}_m$ is a tensor of type (1,1).
The  metrical  linear  connection
${\cal RL}$ introduced by (1.2)
 may  be used to  define the
   covariant derivative  $\nabla$ in ${\cal R}^N$
according to the  conventional  rule:
\be
\nabla_iW^n{}_m=\D{W^n{}_m}{x^i}+L^k{}_i   \D{W^n{}_m}{t^k}
+
L^n{}_{hi}W^h{}_m
-L^h{}_{mi}W^n{}_h,
\ee
which can be written shortly with the help of the operator
$d^{\text{Riem}}_i$ defined by  (1.3), namely
\be
\nabla_iW^n{}_m= d^{\text{Riem}}_iW^n{}_m
+
L^n{}_{hi}W^h{}_m
-L^h{}_{mi}W^n{}_h.
\ee
We have
\be
\nabla_iS=0, \qquad \nabla_i y^j=0, \qquad  \nabla_ia_{mn}=0.
\ee

Due to  the nullifications $\cD_iy^n_h=0$
and
$\cD_i t^j=0$
(see (2.39) and (2.43)),
 we have the transitivity property
\be
\cD_iw^n{}_m=y^n_h t^j_m\nabla_iW^h{}_j
\ee
for the covariant derivatives.

In the commutator
\be
\lf[\nabla_i\nabla_j-\nabla_j\nabla_i\rg] W^n{}_k
=
-y^ma_m{}^h{}_{ij}  \D{W^n{}_k}{y^h}
-
a_k{}^h{}_{ij}W^n{}_h
+
a_h{}^n{}_{ij}W^h{}_k
\ee
the associated Riemannian curvature tensor is constructed in the ordinary way
\be
a_n{}^i{}_{km}=
\D{a^i{}_{nm}}{x^k}-\D{a^i{}_{nk}}{x^m}+a^u{}_{nm}a^i{}_{uk}-a^u{}_{nk}a^i{}_{um}.
\ee

With the ordinary Riemannian covariant derivative
\be
\nabla_k a_h{}^t{}_{ij}
=
\D{a_h{}^t{}_{ij}}{x^k}
+
 a^t{}_{ku}  a_h{}^u{}_{ij}
- a^u{}_{kh}  a_u{}^t{}_{ij}
- a^u{}_{ki}  a_h{}^t{}_{uj}
- a^u{}_{kj}  a_h{}^t{}_{iu},
\ee
the cyclic identity
\be
\nabla_ka_m{}^n{}_{ij}+\nabla_ja_m{}^n{}_{ki}+\nabla_ia_m{}^n{}_{jk}=0
\ee
holds.

Under these conditions,
by comparing the Finslerian commutator (4.13)
 with the  Riemannian precursor
(4.21),
we obtain
\be
M^n{}_{ij}=
-
y^n_t t^h
a_h{}^t{}_{ij}
\ee
\ses
and
\be
E_k{}^n{}_{ij}=
y_h^n t^h_{km}M^m{}_{ij}
+
y^n_ma_h{}^m{}_{ij} t^h_k,
\ee
\ses
together with
$$
\Rho_k{}^n{}_{ij}
=
\lf(
y_h^n t^h_{km}
-C^n{}_{mk}
\rg)
M^m{}_{ij}
+
y^n_ma_h{}^m{}_{ij} t^h_k.
$$
Inserting    here the tensor $C^n{}_{mk}$
taken from (2.22)
and noting
 the vanishing $l_mM^m{}_{ij}=0$
(see (4.6)), we get

$$
\Rho_k{}^n{}_{ij}
=
\lf(
y_h^n t^h_{km}
-
(1-H)\fr1{F}(l_k\de^n_m+l^ng_{mk})
-p^2 t^l_mt^h_{rk} a_{lh}g^{nr}
\rg)
M^m{}_{ij}
+
y^n_ma_h{}^m{}_{ij} t^h_k.
$$

Let us lower here the index $n$
and use the equality
$
g_{nm}y^m_i=p^2t^j_n a_{ij}
$
(ensued from (2.14)).
This yields
$$
\Rho_{knij}
=
\lf(
p^2t^l_n t^h_{km}
a_{lh}
-
(1-H)\fr1{F}(l_kg_{mn}+l_ng_{mk})
-p^2 t^l_mt^h_{nk} a_{lh}
\rg)
M^m{}_{ij}
+
p^2
a_{ml}
t^l_n
a_h{}^m{}_{ij} t^h_k.
$$
Next, we use here the skew-symmetry relation (2.21), obtaining
$$
\Rho_{knij}
=
\lf(
(1-H)\fr2{F}(l_ng_{mk}-l_mg_{kn})
-
(1-H)\fr1{F}(l_kg_{mn}+l_ng_{mk})
\rg)
M^m{}_{ij}
+
p^2
a_{ml}
t^l_n
a_h{}^m{}_{ij} t^h_k,
$$
\ses\\
or
\be
\Rho_{knij}
=
-
(1-H)\fr1{F}(l_kM_{nij}-l_nM_{kij})
+
p^2
a_{hlij} t^h_kt^l_n,
\ee
\ses\\
where
$
a_{hlij}=a_{lr}a_h{}^r{}_{ij}.
$
Finally, we return the index $n$ to the upper position, arriving at
\be
\Rho_k{}^n{}_{ij}
=
-
(1-H)\fr1{F}(l_k\de^n_m-l^ng_{mk})
M^m{}_{ij}
 +
y^n_ma_h{}^m{}_{ij} t^h_k.
\ee

The totally contravariant representation
$$
\Rho^{knij}=g^{pk}a^{mi}a^{nj}
\Rho_p{}^n{}_{mn}
$$
reads
\ses\\
\be
\Rho^{knij}
=
-
(1-H)\fr1{F}(l^kM^{nij}-l^nM^{kij})
 +
\fr1{p^2}
 y^k_h y^n_ra^{hrij},
\ee
where
$a^{hrij}=
a^{hl}a^{mi}a^{nj}
a_{l}{}^r{}_{mn}
$
and
$
M^{mij}=
a^{hi}a^{nj}
M^m{}_{hn}.
$

Similarly, we can conclude  from (4.25) that
the tensor
$$ M_{nij}= g_{nm}M^m{}_{ij}$$
 reads
\be
M_{nij}=-p^2t^ht^m_na_{hmij}.
\ee
Squaring yields
\ses
\be
 M^{nij} M_{nij}=
p^2
t^l
a_l{}^{nij}
t^h
a_{hnij}.
\ee

From the representations (4.27) and (4.30)
it follows directly that the cyclic identities
(4.9) and (4.15) are  consequences
of
the Riemannian
cyclic identity (4.24), for
$\cD_l F=  \cD_l l_k=\cD_l t^h_k=\cD_l p=\cD_l t^m=0$.

Now we square the $\rho$-tensor:
$$
\Rho^{knij}
\Rho_{knij}
=
(1-H)^2\fr2{F^2}M^{nij}M_{nij}
-
2(1-H)\fr1{F}(l^kM^{nij}-l^nM^{kij})
p^2
a_{hlij} t^h_kt^l_n
+
a^{knij}   a_{knij}
$$

\ses

\ses

\ses

$$
=
(1-H)^2\fr2{F^2}M^{nij}M_{nij}
-
2(1-H)H\fr1{F^2}
p^2
(a_{hlij} t^ht^l_nM^{nij}-a_{hlij} t^h_kt^lM^{kij})
+
a^{knij}   a_{knij},
$$
\ses
or
$$
\Rho^{knij}
\Rho_{knij}
=
(1-H)^2\fr{2p^2}
{F^2}
t^l
a_l{}^{nij}
t^h
a_{hnij}
+
2(1-H)\fr{Hp^2}{F^2}
(a_{hlij} t^ht^ra_r{}^{lij}-a_{hlij} t^lt^ra_r{}^{hij})
+
a^{knij}   a_{knij},
$$
\ses
which is
\be
\Rho^{knij}\Rho_{knij}=
a^{knij}a_{knij}
+\fr 2{S^2}
\lf(\fr1{H^2}-1\rg)
t^l
a_l{}^{nij}
t^h
a_{hnij}.
\ee

Because of the transitivity  (4.20),
from (4.25) it follows that
\be
\cD_l M^n{}_{ij}
=
-
y^n_tt^h
\nabla_l a_h{}^t{}_{ij}.
\ee

From (4.28) we can conclude that
\be
\cD_l \Rho_k{}^n{}_{ij}=
(1-H)\fr1{F}(l_k\de^n_m-l^ng_{mk})
y^m_tt^h
\nabla_l a_h{}^t{}_{ij}
 +
y^n_mt^h_k\nabla_la_h{}^m{}_{ij} .
\ee

It is also convenient to use the representation
\be
\Rho_{knij}=T_{kn}{}^{hm}a_{hmij},
\ee
where
\be
T_{kn}{}^{hm}= p^2
\lf[
\fr12(t^h_kt^m_n-t^m_kt^h_n)
+(1-H)\fr1{F^2}(y_kt^ht^m_n-y_nt^ht^m_k)
\rg].
\ee
\ses
Since
\be
\cD_lT_{kn}{}^{hm}=0,
\ee
we have the relation
\be
\cD_l\Rho_{knij}=T_{kn}{}^{hm}\nabla_la_{hmij}.
\ee

 \ses

\ses

\setcounter{sctn}{5} \setcounter{equation}{0}

\ses

\ses

\bc   {\bf  5. ${\cal FS}$-space  example of the space ${\cal F}^N$ } \ec

\ses

\ses

Let us also assume that the manifold $M$ admits a non--vanishing 1-form $b= b(x,y)$
and
denote by
\be
c=|| b||_{\text{Riemannian}}
\ee
the respective Riemannian norm value,
assuming
\be
 0 < c < 1.
 \ee

\ses

With respect to  natural local coordinates $x^i$
we have the local representations
\be
a^{ij}(x)b_i(x)b_j(x)=c^2(x),  \qquad
 b=b_i(x)y^i.
\ee
The reciprocity  $a^{in}a_{nj}=\de^i{}_j$ is assumed, where $\de^i{}_j$ stands for the Kronecker symbol.
The covariant index of the vector $b_i$  will be raised by means of the Riemannian rule
$ b^i=a^{ij}b_j,$ which inverse reads $ b_i=a_{ij}b^j.$

\ses

We shall use also the normalized vectors
\be
\wt b_i=\fr1cb_i, \qquad  \wt b^i=\fr1cb^i= a^{ij} \wt b_j, \qquad
 a^{mn}  \wt b_m \wt b_n=1.
\ee

We get
\be
a_{ij}y^iy^j-b^2>0
\ee
and may conveniently use the variable
\be
q:=\sqrt{a_{ij}y^iy^j-b^2}.
\ee
Obviously, the inequality
\be
q^2 \ge \fr{1-c^2}{c^2}\,b^2
\ee
is valid.

We
also  introduce the tensor
\be
r_{ij}(x)~:=a_{ij}(x)-b_i(x)b_j(x)
\ee
to have the representation
\be
q=\sqrt{r_{ij}y^iy^j}.
\ee
The equalities
\be
r_{ij}b^j=(1-c^2)b_i, \qquad
  r_{in}r^{nj}=r^j{}_i-(1-c^2)b^jb_i
 \ee
hold.

In evaluations it is convenient to use the variables
\be
u_i~:=a_{ij}y^j,
\qquad
v^i~:=y^i-bb^i, \qquad v_m~:=u_m-bb_m=r_{mn}y^n\equiv a_{mn}v^n.
\ee
We have
\be
r_{ij}=\D{v_i}{y^j},
\qquad
\D b{y^i}=b_i, \qquad \D q{y^i}=\fr{v_i}q,
 \qquad
v_ib^i=v^ib_i=(1-c^2)b,
\ee
\ses
and
\be
u_iv^i=v_iy^i=q^2,
\qquad
 r_{in}v^n=v_i-(1-c^2)bb_i,  \qquad
v_kv^k=q^2-(1-c^2)b^2.
\ee

With  the variable
\be
w=\fr qb,\qquad b>0,
\ee
we obtain
\be
\D{w}{y^i}=\fr{qe_i}{b^2}, \qquad e_i=-b_i+\fr b{q^2}v_i,
\ee
and
$
y^ie_i=0.
$

The
Finsler metric function $F$ of the ${\cal FS}$-space
is given by (1.25).
When $b>0$, we can conveniently
use  the {\it generating metric function} $V=V(x,w)$ to have the representation
\be
F=bV(x,w).
\ee

The unit vector
$l_m=\partial F/\partial y^m$
is given by
\be
l_m=b_mV+we_mV', \qquad
V'= \D Vw.
\ee

It proves convenient to use the quantities
\be
\tau= \fr{wV}{V'},
\ee
\ses
\be
\wt  q= \sqrt{q^2+\lf(1-\fr1{c^2}\rg)b^2},
\ee
\ses
and
\be
  \qquad \wt w= \fr{\wt q}b,   \qquad  b>0.
\ee

There are  the useful equalities
\ses
$$
 \tau= \wt \tau =\fr{\wt wV}{\wt V'},
\qquad \wt V= V,
\qquad
\wt V'= \D {\wt V}{\wt w},
\qquad
b^ml_m=c^2V\lf(1-\fr{\wt w^2}{\tau}\rg).
$$

\ses

We say that the ${\cal FS}$-space is {\it special}, if
$\partial\wt  V/\partial x^n=0$, that is when
\be
\wt V=\wt V(\wt w).
\ee

Take two differentiable scalar functions
\be
C=C(x), \quad C_1= C_1(x)     , \qquad  C>0, \quad C>|C_1|,
\ee
\ses
and construct the scalars
\be
H=\sqrt{C^2-(C_1)^2}
\ee
and
\be
 \breve k
=
\sqrt{\fr{C-C_1}{C+C_1}}.
\ee
\ses
Let a positive function $\mu=\mu(x,y)$
be specified according to
\be
\sqrt{\mu}
=
\fr H{2 \breve k}
\bigl[
1+\breve k^2 +(1-\breve k^2  )   \cos\varrho
\bigr],
 \ee
where
 $\varrho=\varrho(x,y)$ is an input scalar.
We can write
\be
\sqrt{\mu}
=
C+C_1\cos\varrho.
\ee

Consider the transformation
 $t^m=t^m(x,y)$
with
\be
  t^m
=
\lf[
i^m \sin\varrho
+
\fr 1{2 \breve k}
[
1-\breve k^2 +(1+\breve k^2  )   \cos\varrho]
\wt b^m
\rg]
\fr{H}{\sqrt{\mu} }
F^H,
\ee
where
\ses\\
\be
i^m=\lf(y^m-\fr1{c^2}
bb^m\rg)
\fr 1{\wt q}.
\ee

\ses

We have
\be
b_mi^m=0, \qquad a_{mn} i^m   i^n  =1,   \qquad a_{mn} y^m   i^n  = \wt q,
\ee
and
\be
 b^*=
\fr 1{2 \breve k}
 \Bigl[
1-\breve k^2
+ (1+\breve k^2)  \cos\varrho
\Bigr]
\fr{H}{\sqrt{\mu} }
S,
\ee
\ses
where
\be
S=\sqrt{a_{mn}t^mt^n}, \qquad     b^*=t^m\breve b_m, \qquad \breve b_m=\wt b_m,
\qquad \breve b^m=\wt b^m.
\ee
We  get  also the equality
\be
 b^*=
(C_1+C\cos\varrho)
\fr{1}{\sqrt{\mu} }
S.
\ee

The functions (5.27) obviously fulfill the $H$-degree homogeneity
condition (2.7).
The validity of the equality
$
S=F^H
$
(see (2.10))
can readily be verified.

\ses

The property
\be
t^m(x,b(x)) \thicksim b^m(x)
\ee
holds.

The following useful equalities can readily be obtained:
$$
 \cos\varrho
 =
 -
\fr
{ (1-\breve k^2)S-(1+\breve k^2) b^* }
{ (1+\breve k^2)S-(1-\breve k^2) b^* },
$$

\ses

\ses

\be
\sqrt{\mu}= \fr{2H\breve k S}{(1+\breve k^2) S-(1-\breve k^2) b^*},
\ee

\ses

\ses

$$
 \cos\varrho
 =
 -
\fr{\sqrt{\mu} }
{2H\breve k S}
[ (1-\breve k^2)S-(1+\breve k^2) b^* ],
$$
\ses
and
$$
 \sin^2\varrho
 =
4\breve k^2
\fr
{ S^2-(b^*)^2 }
{\bigl[ (1+\breve k^2)S-(1-\breve k^2) b^* \bigr]^2},
$$
\ses\\
together with
\be
\fr{ \sin^2\varrho}{\mu}=
\fr{1}{H^2}
\fr
{ S^2-(b^*)^2 }
{S^2}.
\ee

Differentiating (5.28) yields
\be
\D{i^m}{y^k}=
\lf(\de^m_k-\fr1{c^2}  b_kb^m\rg)
\fr 1{b\wt w}
-\fr1bb_ki^m
-
\fr 1b
\fr { w^2}{\wt w^2}
i^m
e_k,
\ee
which entails
\ses\\
\be
y^k\D{i^m}{y^k}=0,  \qquad   b^k\D{i^m}{y^k}=0,
\qquad
i^na_{nm}\D{i^m}{y^k}=0,  \qquad   b_m\D{i^m}{y^k}=0,
\qquad
y^na_{nm}\D{i^m}{y^k}=0.
\ee

\ses

We can use the relations
$$
b_k=\fr bFl_k-w^2\fr1{\tau}e_k, \qquad
\wt w  i_k
=
w^2e_k + \wt w^2b_k
=
w^2e_k + \wt w^2
 \lf(\fr bFl_k-w^2\fr1{\tau}e_k\rg),
 $$
so that
$$
i_k
-
\fr1F b\wt w
l_k
=
\fr{w^2}{\wt w}
\lf(1-\fr{\wt w^2}{\tau}\rg)
e_k,
 $$
where
$
i_k=a_{kn}i^n.
$

We have also
$$
a_{mh}\D{i^m}{y^k}=
\lf(a_{hk}-\fr1{c^2}  b_kb_h\rg)
\fr 1{b\wt w}
-\fr1bb_ki_h
-
\fr 1b
\fr { w^2}{\wt w^2}
i_h
e_k,
$$
\ses
which entails
\be
a_{mh}\D{i^m}{y^k}=
\lf(a_{hk}-\fr1{c^2}  b_kb_h-i_ki_h\rg)
\fr 1{b\wt w}.
\ee

With these observations, from (5.27) we find that the
derivative coefficients
$    t^m_k=\partial{t^m}/\partial{y^k}   $
can be given by
$$
\fr1H
t^m_k=
\lf[\cos\varrho i^m
-
\fr 1{2 \breve k}
(1+\breve k^2  )
  \sin\varrho
\breve b^m
\rg]
\varrho'
\fr wbe_k
\fr{F^H}{\sqrt{\mu} }
+
\sin\varrho
\D{i^m}{y^k}
\fr{F^H}{\sqrt{\mu} }
+
\fr 1Fl_k
t^m
-
\fr{1}{H\sqrt{\mu} }
\D{\sqrt{\mu}}{y^k}
t^m.
$$
\ses
Since
$$
\D{\sqrt{\mu}}{y^k}=
-
\fr H{2 \breve k}
(1-\breve k^2  )
 \sin\varrho\,
\varrho'
\fr wbe_k,
$$
 we obtain the explicit representation
\ses\\
$$
\fr1H
\sqrt{\mu}
t^m_k=
\lf[\cos\varrho i^m
-
\fr 1{2 \breve k}
(1+\breve k^2  )
  \sin\varrho
\breve b^m
\rg]
\varrho'
\fr wbe_k
F^H
$$

\ses

\ses

\be
+
\sin\varrho
\D{i^m}{y^k}
F^H
+
\lf[
 \sqrt{\mu}
 \fr 1Fl_k
+
\fr 1{2 \breve k}
(1-\breve k^2  )
 \sin\varrho\,
\varrho'
\fr wbe_k
\rg]
t^m.
\ee
\ses\\
The identity
$
t^m_ky^k=Ht^m
$
can readily be verified.

We can straightforwardly
evaluate the contraction
$a_{mn} t^m_k   t^n_h $,
which leads to the expression which is a linear combination of
$g_{kh}$, $e_ke_h$,  $l_kl_h$, and
$e_kl_h+e_hl_k$.
To obtain the conformal result, the terms $l_kl_h$ are to be canceled,
which proves possible if and only if the function $\mu$ is taken to be
\ses\\
\be
\mu=\fr{1}{\wt w^2}\tau\sin^2\varrho,
\ee
which entails
\ses\\
\be
\fr{\wt w^2}{\tau}=
\fr{1}{H^2}
\fr
{ S^2-(b^*)^2 }
{S^2}
\ee
(see (5.35)).
With the choice of $\mu$ according to  (5.40), using the representation (5.39)
leads  straightforwardly to the equality
$$
\fr{1}{H^2}
a_{mn} t^m_k   t^n_h  =
\fr1{\tau}
\fr{\wt w^2 }{\sin^2\varrho}
\lf( \D{\varrho}{ w}\rg)^2
\fr{w^2}{b^2} e_ke_h
F^{2H}
 $$

\ses

\ses

\ses

\be
+
\fr1{\tau}
 F^{2H}
\lf(
-
 \fr{\tau- w(\tau'- w)}   {\tau}
+1-\fr{w^2}{\wt w^2}
\rg)
w^2e_ke_h
\fr1{b^2}
+
   F^{2(H-1)}
g_{kh}
\ee
(see  Appendix E).

\ses

Subjecting the $\varrho$ to the equation
\be
 \D{\varrho}{\wt w}=
\fr 1{\wt w}
\sqrt
{\fr{\wt\tau-\wt w(\wt\tau'-\wt w)}
{\wt\tau}}
\sin\varrho,
\ee
where $\wt\tau'=\partial \tau/\partial \wt w$,
is necessary and sufficient to reduce the right-hand part in  (5.42)
to the conformal representation, namely we obtain
simply
\be
\fr{1}{H^2}
a_{mn} t^m_k   t^n_h  =
   F^{2(H-1)}
g_{kh}.
\ee

Comparing the last representation with
 the formulas (2.1), (2.4), and (2.9)
makes us conclude that  the  following assertion is valid.

\ses

\ses

{\bf Proposition 5.1.}
{\it With choosing the function $\mu$ to be given by} (5.40)
{\it and subjecting the function
$\varrho$ to the equation}
 (5.43),
{\it the transformation}
 $t^m=t^m(x,y)$
 {\it introduced  by}  (5.27)
{\it fulfills the conformal automorphism} (2.1)-(2.2).

\ses

\ses

From (5.28) it follows that
$$
\fr 1{b}
y^m
=
\fr{\sqrt{\tau} }
{H S}
(  t^m-b^*\breve b^m)
+ \fr1{c^2}b^m,
$$
so that we
can write (5.27) to read
\ses\\
$$
  t^m
-
b^*\breve b^m
=
\fr{H F^H}
{\sqrt{\tau} }
\lf[
\fr 1{b}y^m
-
 \fr1{c^2}b^m
\rg].
$$
On the other hand,
using the equality $S=F^H$ in (5.41) yields
\ses\\
\be
S^2-(b^*)^2=
\fr{H^2 F^{2H}}
{\tau }
\wt w^2.
\ee
\ses
Thus,  it is valid that
\ses\\
\be
\fr1{\wt w}
\lf(
\fr 1{b}y^m
-
 \fr1{c^2}b^m
\rg)
=
\fr1{\sqrt{S^2-(b^*)^2}}
\lf(  t^m-b^*\breve b^m\rg).
\ee

Let us evaluate the explicit form of the respective coefficients
$N^m{}_n$ proposed by (2.36).
Accordingly, we assume
  $d_iF=0$   and $H=const$.
Since $S=F^H$, we have
$
d_iS=0.
$
 From (2.38)--(2.39) it follows that
$
d_it^n=-a^n{}_{ih}t^h.
$

Also,
$
d_ib^*=d_i(t^m\wt b_m)= -a^m{}_{ih}t^h\wt b_m+t^m(\nabla_i\wt b_m+a^r{}_{im}\wt b_r)
=t^m\nabla_i\wt b_m,
$
\ses
so that
$$
\fr1{\sqrt{S^2-(b^*)^2}}
 d_ib^*
=
\lf[
\fr1{\wt w}
\lf(
\fr 1{b}y^m
-
 \fr1{c^2}b^m
\rg)
+
\fr1{\sqrt{S^2-(b^*)^2}}
b^*\breve b^m
\rg]
\nabla_i\wt b_m.
$$
\ses
Henceforth, we assume $c=const$. In this case
$
b^m\nabla_ib_m=0,
$
and therefore
$$
\fr1{\sqrt{S^2-(b^*)^2}}
 d_ib^*
=
\fr1{\wt w}
\fr 1{b}y^m
\nabla_i\wt b_m.
$$

\ses

Under these conditions, applying the operator $d_i$
to (5.46)
yields
\ses\\
$$
-
\fr1{\wt w^2}
d_i\wt w
\lf(
\fr 1{b}y^m
-
 \fr1{c^2}b^m
\rg)
+
\fr1{\wt w}
\lf(
\fr 1{b}
d_iy^m
-
\fr 1{b^2}
y^md_ib
-
 \fr1{c^2}
( \nabla_ib^m-a^m{}_{ih}b^h)
\rg)
$$

\ses

$$
=
b^*d_ib^*
\fr1{(S^2-(b^*)^2)\,\sqrt{S^2-(b^*)^2}}
\lf(  t^m-b^*\breve b^m\rg)
$$

\ses

\ses

$$
+
\fr1{\sqrt{S^2-(b^*)^2}}
\lf( d_i t^m- (d_ib^*)\breve b^m
-
b^*
(\nabla_i\breve b^m    -a^m{}_{ih}\wt b^h)
\rg),
$$
\ses

\nin
or
\ses\\
$$
-
\lf[
\fr1{\wt w}
d_i\wt w
+
\fr1
{S^2-(b^*)^2}
b^*d_ib^*
\rg]
\lf(
\fr 1{b}y^m
-
 \fr1{c^2}b^m
\rg)
+
\fr 1{b}
d_iy^m
-
\fr 1{b^2}
y^md_ib
-
 \fr1{c^2}
( \nabla_ib^m-a^m{}_{ih}b^h)
$$

\ses

$$
=
\fr{\wt w}{\sqrt{S^2-(b^*)^2}}
\lf(
 -a^m{}_{ih}
 \lf(
b^*\breve b^h
+
\fr{H F^H}
{\sqrt{\tau} }
\lf(
\fr 1{b}y^h
-
 \fr1{c^2}b^h
\rg)
\rg)
 -b^*(\nabla_i\breve b^m    -a^m{}_{ih}\wt b^h)
\rg)
-
\fr1{b}
(y^h\nabla_i\breve b_h)
 \breve b^m.
$$

Using the new variable
\be
W=\fr{b^*}{S},
\ee
\ses
we have
\be
 \cos\varrho
 =
 -
\fr
{ 1-\breve k^2-(1+\breve k^2) W }
{ 1+\breve k^2-(1-\breve k^2) W }
\ee
and
\ses\\
\be
 \sin^2\varrho
 =
4\breve k^2
\fr
{ 1-W^2 }
{\bigl[ 1+\breve k^2-(1-\breve k^2) W \bigr]^2},
\ee
which entails
\ses\\
\be
\D{ \cos\varrho}W
 =
4\breve k^2
\fr
1
{\bigl[ 1+\breve k^2-(1-\breve k^2) W \bigr]^2}.
\ee

\ses

In the equality
$$
\fr1{\wt w}\D {\wt w}W
=
\fr1{\sin\varrho}
\D{\varrho}W
P
$$
we use (5.49) and (5.50), obtaining
\ses
$$
\fr1{\wt w}\D {\wt w}W
=
-
\fr1{1-W^2}
P,
$$
where the notation
$$
P=\sqrt
{\fr{\wt\tau}{\wt\tau-\wt w(\wt\tau'-\wt w)}
}
$$
has been used.

Therefore,
in the special case
(see (5.21))
of the ${\cal FS}$-space, we have
$$
d_i\wt w=
\D {\wt w}W
d_iW
=
-
\fr{\wt w}{1-W^2}
P
d_iW
=
-
\fr{\wt w}{1-W^2}
P
\fr1S
d_ib^*
$$
\ses
and
can write
$$
\fr1S(d_ib^*)
\fr1{1-W^2}
(P-W)
\lf(
\fr 1{b}y^m
-
 \fr1{c^2}b^m
\rg)
+
\fr 1{b}
d_iy^m
-
\fr 1{b^2}
y^md_ib
-
 \fr1{c^2}
 \nabla_ib^m
$$

\ses

$$
=
\fr{\wt w}{\sqrt{S^2-(b^*)^2}}
\lf(
 -a^m{}_{ih}
\fr{H F^H}
{\sqrt{\tau} }
\fr 1{b}y^h
 -b^*\nabla_i\breve b^m
\rg)
-
\fr1{b}
(y^h\nabla_i\breve b_h)
 \breve b^m
$$

\nin
and
$$
d_iy^m=-
\fr1S(d_ib^*)
\fr1{1-W^2}
(P-W)
\lf(
y^m
-
 \fr1{c^2}bb^m
\rg)
+
\fr 1{b}
y^md_ib
+
 \fr1{c^2}
b \nabla_ib^m
$$

\ses

\ses

$$
+
\fr{\wt w}{\sqrt{S^2-(b^*)^2}}
\lf(
 -a^m{}_{ih}
\fr{H F^H}
{\sqrt{\tau} }
y^h
 -bb^*\nabla_i\breve b^m
\rg)
-
(y^h\nabla_i\breve b_h)
 \breve b^m.
$$

Taking into account the equality
$$
d_ib=-\fr{b\wt w}{\tau}d_i\wt w
$$
which is valid
in the special case of the ${\cal FS}$-space,
we obtain
\ses\\
$$
d_iy^m=-
\fr1S(d_ib^*)
\fr1{1-W^2}
(P-W)
\lf(
y^m
-
 \fr1{c^2}bb^m
\rg)
+
\fr{\wt w^2}{\tau}
\fr1S(d_ib^*)
\fr1{1-W^2}
Py^m
$$

\ses

\ses

$$
+
\fr{\sqrt{\tau }} {H F^{H}}
\lf(
 -a^m{}_{ih}
\fr{H F^H}
{\sqrt{\tau} }
y^h
 -bb^*\nabla_i\breve b^m
\rg)
-
(y^h\nabla_i\breve b_h)
 \breve b^m
+
 \fr1{c^2}
b \nabla_ib^m
$$

\ses

\ses

\ses

$$
=
-
            \fr1{b \wt w}         y^h\nabla_i\breve b_h
\fr1{\sqrt{1-W^2}}
(P-W)
\lf(
y^m
-
 \fr1{c^2}bb^m
\rg)
$$

\ses

\ses

$$
+
\fr{\wt w^2}{\tau}
            \fr1{b \wt w}         y^h\nabla_i\breve b_h
\fr1{\sqrt{1-W^2}}
Py^m
-
b
\fr{\sqrt{\tau }} {H }
W
\nabla_i\breve b^m
-
(y^h\nabla_i\breve b_h)
 \breve b^m
+
 \fr1{c^2}
b \nabla_ib^m
 -
 a^m{}_{ih}
y^h.
$$
\ses
Here,
$$
1-W^2=H^2\fr{\wt w^2}{\tau}.
$$

Therefore,
\ses\\
$$
d_iy^m=-
            \fr1{b \wt w}         y^h\nabla_i\breve b_h
\fr1H\fr{\sqrt{\tau}}{\wt w}
(P-W)
\lf(
y^m
-
 \fr1{c^2}bb^m
\rg)
+
\fr{\wt w^2}{\tau}
            \fr1{b \wt w}         y^h\nabla_i\breve b_h
\fr1H\fr{\sqrt{\tau}}{\wt w}
Py^m
$$

\ses

\ses

\be
-
b
\fr{\sqrt{\tau }} {H }
W
\nabla_i\breve b^m
-
(y^h\nabla_i\breve b_h)
 \breve b^m
+
 \fr1{c^2}
b \nabla_ib^m
 -
 a^m{}_{ih}
y^h.
\ee

\ses

Noting the equality
\be
N^m{}_i=d_iy^m
\ee

\nin
leads to
\ses\\
\be
N^m{}_i
=
\fr1H
\fr1{b\wt w}      (   y^h\nabla_i\wt b_h)
F\al^m
-
\fr1H
 \sqrt{\tau-H^2\wt w^2 }    \,
b\wt\beta^m_i
-
(y^h\nabla_i\wt b_h)
 \wt b^m
+
\wt b \nabla_i \wt b^m
 -
 a^m{}_{ih}
y^h,
\ee
\ses
\ses\\
 where
\be
F\al^m
=
\fr{\wt w}{\sqrt{\wt\tau-\wt w(\wt\tau'-\wt w)}
}
\lf[
y^m
-
\fr{\tau}{\wt w^2}
\lf(
y^m
-
\wt b\wt b^m
\rg)
\rg]
\ee
 and
\ses\\
\be
\wt\beta^m_i
=
\nabla_i\wt b^m
-
            \fr1{b^2 \wt w^2}  (       y^h\nabla_i\wt b_h)
\lf(
y^m
-
\wt b\wt b^m
\rg),
\ee

\nin
which  can also be written in the form
\ses\\
\be
N^m{}_i
=
\fr1H
\fr1{\wt q}      (   y^h\nabla_i \wt b_h)
F\al^m
-
\fr1H
 \sqrt{B-H^2\wt q^2 }    \,
\wt\beta^m_i
-
(y^h\nabla_i \wt b_h)
\wt b^m
+
\wt b \nabla_i \wt b^m
 -
 a^m{}_{ih}
y^h,
\ee
\ses
\ses\\
with
\be
F\al^m
=
\fr{\wt w}{\sqrt{\wt\tau-\wt w(\wt\tau'-\wt w)}
}
\lf[
y^m
-
\fr{B}{\wt q^2}
\lf(
y^m
-
\wt b \wt b^m
\rg)
\rg],
\ee
\ses\\
\be
\wt\beta^m_i
=
\nabla_i \wt b^m
-
            \fr1{ \wt q^2}  (       y^h\nabla_i\wt b_h)
\lf(
y^m
-
\wt b \wt b^m
\rg),
\ee
\ses\\
and
\be
B=b^2\tau.
\ee

Thus we have

\ses

\ses

{\bf Proposition 5.2.}
{\it If in the special case of the ${\cal FS}$-space with $c=const$
the transformation} (5.27) {\it results in the conformally automorphic space,
then  the coefficients} (2.36) {\it can explicitly be given by means of
the representation}
(5.56)-(5.59).

\ses

\ses

It is easy to verify that
\ses
\be
\al^ml_m=0,  \qquad
\wt\beta^m_i b_m=0,  \qquad  \wt\beta^m_il_m=0, \qquad  \wt\beta^m_i\al_m=0.
\ee

By contracting (5.56) we find
\be
 l_mN^m{}_i=
-
(   y^h\nabla_i b_h)
\lf(1-\fr{1+w^2}{\tau}
\rg)
V
 -
l_m a^m{}_{ih}
y^h.
 \ee
From this result it follows that
\be
\D {F}{x^k}+ l_mN^m{}_k=0.
\ee

Indeed,
denoting
\be
s_k~:      =y^m\nabla_kb_m,
\ee
we get
\be
\D {q}{x^k}=-\fr b{q}(s_k+y^mb_ha^h{}_{mk})
+\fr1qy^my^n     \D{a_{mn}}{x^k}
\ee
and
$$
\D {w}{x^k}=
-\fr 1{q}(s_k+y^mb_ha^h{}_{mk})
-\fr q{b^2}(s_k+y^mb_ha^h{}_{mk})
+\fr1{bq}  y^my^n  \D{a_{mn}}{x^k},
$$
or
\be
\D {w}{x^k}=
-\fr 1{b^2q}S^2(s_k+y^mb_ha^h{}_{mk})
+\fr1{bq}  y^my^n  \D{a_{mn}}{x^k}.
\ee
The equality
\be
\wt w \D {\wt w}{x^k}=w \D {w}{x^k}
\ee
can appropriately be used.

\ses

In the  special  case
$
F=b \wt V(\wt w)
$
(see (5.21))
of the ${\cal FS}$-space
we have
\ses\\
$$
\D {F}{x^k}=    \lf(\wt V   -\fr 1{bq}S^2 \fr w{\wt w}    \wt V'\rg)s_k
+\lf(\wt V   -\fr 1{bq}S^2 \fr w{\wt w}    \wt V'\rg)y^mb_ha^h{}_{mk}
+\fr1{q} \fr w{\wt w}     \wt V' y^my^n  \D{a_{mn}}{x^k},
$$
\ses\ses
or
\be
\D {F}{x^k}=    \lf(\wt V   -\fr {S^2}{b^2} \fr 1{\wt w}    \wt V'\rg)s_k
+\lf(\wt V   -\fr {S^2}{b^2} \fr 1{\wt w}    \wt V'\rg)y^mb_ha^h{}_{mk}
+\fr1{b\wt w}     \wt V' y^my^n  \D{a_{mn}}{x^k}.
\ee

\ses

\ses

In terms of the function $\wt\tau$, we can write
\ses
\be
\D {F}{x^k}=  \wt V  \lf(1   -\fr {1+w^2}{\wt\tau}\rg)  s_k
+\wt V  \lf(1   -\fr {1+w^2}{\wt\tau}\rg)  y^mb_ha^h{}_{mk}
+ \fr {\wt V}{b\wt\tau} y^my^n  \D{a_{mn}}{x^k}.
\ee
With  this equality the validity of the vanishing (5.62) can readily be verified.

The following proposition  is valid.

\ses

\ses

{\bf Proposition 5.3.}
{\it The  transformation} (5.27)
{\it entails the conformal automorphism}
(2.1)
{\it  iff}
\be
\tau={\breve C}^2+2{\breve C}\sqrt{1-H^2}\,\wt w+\wt w^2.
\ee

\ses

It follows that
$$
\wt\tau-\wt w(\wt\tau'-\wt w)={\breve C}^2.
$$

\ses

In these formulas, ${\breve C}$ is an integration scalar ${\breve C}={\breve C}(x)$.
It can readily be seen that when $|{\breve C}|\ne 1$,
the entailed Finsler metric function $F$
can vanish at various values of tangent vectors $y$.
To agree with the condition that $F$ vanishes only at zero-vectors $y=0$,
we admit  strictly the values ${\breve C}=1$ and ${\breve C}=-1$.
In this case we can write the above $\tau$ as follows:
\be
\tau=1+g\wt w+\wt w^2, \quad -2<g<2.
\ee
Generally, the $g$ may depend on $x$.
We obtain
\be
B-H^2\wt q^2=
\lf(b+\fr12g\wt q\rg)^2.
\ee

In this case the coefficients (5.56) take on the form
\ses\\
\be
N^m{}_i
=
\fr1h
\fr1{\wt q}
Fm^m
\wt s_i
-
\fr1h
\lf(b+\fr12g\wt q\rg)
\wt\beta^m_i
-
\wt b^m
\wt s_i
+
\wt b \nabla_i \wt b^m
 -
 a^m{}_{ih}
y^h,
\ee
\ses
\ses\\
with
\be
m^m
=
\fr1{ \wt q F}
\lf[
 \wt q^2 \wt b^m
-
\lf(b+g\wt q\rg)
\lf(
y^m
-
\wt b \wt b^m
\rg)
\rg] \equiv
\fr1{ \wt q F}
\lf[
B^2 \wt b^m
-
\lf(b+g\wt q\rg)
y^m
\rg],
\ee
\ses\\
\be
\wt\beta^m_i
=
\nabla_i \wt b^m
-
            \fr1{ \wt q^2}
\lf(
y^m
-
\wt b \wt b^m
\rg)\wt s_i,
\ee
\ses\\
$m^m=\al^m$,
and
\be
\wt s_i=   y^h\nabla_i \wt b_h.
\ee

\ses

\ses

{\bf Note.}
We  used the input representation $F=bV(x,w)$ (see (5.16)) at $b>0$.
All the performed calculations can be repeated word-for-word in the negative case
$b<0$.
The above representation (5.72)-(5.75) obtained for the coefficients
$N^m{}_i$ embraces both the cases $b>0$ and $b<0$.

\ses

\ses

The last three terms in (5.72) are linear with respect to the tangent vectors $y$.

The function $\tau$ given by (5.70)
 represents the ${\mathbf\cF\cF^{PD}_{g}}$-Finsleroid space
 described in the paper
 [7].
To comply
with the representations used in [7],
we should replace the notation $H$ by the notation $h$:
\be
h=\sqrt{1- \fr{g^2}4}.
\ee
The $g$ plays the role of the characteristic parameter.
The ${\mathbf\cF\cF^{PD}_{g}}$-Finsleroid metric function
$K$
is given as it follows:
\be
K=\sqrt B\, J, \qquad \text{with} ~~ J=\e^{-\frac12 g \chi},
\ee
\ses
where
\be
\chi=\fr1h \Bigl(
-\arctan   \fr G2   +\arctan\fr{L}{hb}\Bigr),    ~  {\rm if}  ~ b\ge 0;
\quad
\chi=\fr1h \Bigl(
 \pi-\arctan
\fr G2
+\arctan\fr{L}{hb}\Bigr),
~   {\rm if}
~ b\le 0,
\ee
 with the function
$    L =\wt q+(g/2) b $
fulfilling   the identity
\be
 L^2+h^2b^2=B.
 \ee
$B$ is the function given by (5.71):
\be
B=b^2+gb\wt q+\wt q^2;
\ee
$G=g/h$.
The definition range
$$
0\le\chi\le\fr1h\pi
$$
is of value to describe all the tangent space.
The normalization in (5.78)
is such that
\be
\chi\bigl|_{y=b}\bigr. =0.
\ee
The quantity $\chi$ can conveniently be written as
\be
\chi  =  \fr1h  f
\ee
with
the function
\be
f=\arccos \fr{ A(x,y)}   {\sqrt{B(x,y)}}, \qquad A=b+\fr12g\wt q,
\ee
ranging as follows:
\be
0\le f\le \pi.
\ee

The function $K$ is the solution for the equation (5.70).

The Finsleroid-axis vector $b^i$ relates to the value $f=0$, and
the opposed vector $-b^i$ relates to the value $f=\pi$:
\be
f=0 ~~  \sim ~~ y=b;  \qquad f=\pi ~~ \sim ~~ y=-b.
\ee
The normalization is such that
\be
K(x,b(x))=1
\ee
(notice that
 $\wt q=0$ at $y^i=b^i$).
The positive  (not absolute) homogeneity  holds:
 $K(x, \ga y)=\ga K(x,y)$ for any $\ga >0$   and all admissible $(x,y)$.

The entailed components $y_i :=(1/2)\partial {K^2}/ \partial{y^i}$
of the  covariant tangent vector $\hat y=\{y_i\}$
can be found in the simple form
\be
y_i=(u_i+ g\wt qb_i) J^2,
\ee
where $u_i=a_{ij}y^j$.

 Under these conditions, we
obtain  the ${\mathbf\cF\cF^{PD}_{g}}$-Finsleroid space
\be
{\mathbf\cF\cF^{PD}_{g}} :=\{M;\,a_{ij}(x);\,b_i(x);\,g(x);\,K(x,y)\}.
\ee

\ses

\ses

 {\large  Definition}.  Within  any tangent space $T_xM$, the  metric function $K(x,y)$
  produces the {\it    ${\mathbf\cF\cF^{PD}_{g} } $-Finsleroid}
 \be
 \cF\cF^{PD}_{g;\,\{x\}}:=\{y\in   \cF\cF^{PD}_{g; \, \{x\}}: y\in T_xM , K(x,y)\le 1\}.
  \ee

\ses

 \ses

 {\large  Definition}. The {\it    ${\mathbf\cF\cF^{PD}_{g} } $-Indicatrix}
 $ {\cal I}\cF^{PD}_{g; \, \{x\}} \subset T_xM$ is the boundary of the
    ${\mathbf\cF\cF^{PD}_{g} } $-Finsleroid, that is,
 \be
{\cal I}\cF^{PD}_{g\, \{x\}} :=\{y\in {\cal I}\cF^{PD}_{g\, \{x\}} : y\in T_xM, K(x,y)=1\}.
  \ee

\ses

 \ses

 {\large  Definition}. The scalar $g(x)$ is called
the {\it Finsleroid charge}.
The 1-form $b=b_i(x)y^i$ is called the  {\it Finsleroid--axis}  1-{\it form}.

\ses

\ses

It can readily be seen that
$$
\det(g_{ij})=\biggl(\fr{K^2}B\biggr)^N\det(a_{ij})>0, \qquad
A^iA_i=     \fr{N^2g^2}4,$$
where $A_i=KC_i.$

\ses

\ses

{\bf Note.}
 The  representation (5.72)-(5.75) obtained for the coefficients
$N^m{}_i$
coincides exactly with the representation (6.53) of [7].
Considering the vector
$C_i=g^{mn}C_{imn}$, the equality
\be
m^m=\fr{C^m}{\sqrt{g^{kh}C_kC_h}}
\ee
holds exactly with the vector $m^m$ given by the representation (5.73)
(which is equivalent to the representation (A.46) proposed in [7]).

\ses

\ses

Let us verify Proposition 5.3.
With the variable
$ W=b^*/S $
(see (5.47)) we can write the equation (5.43) as follows:
\be
\fr1{\sin\varrho}
\D{\varrho}W
\D W{\wt w}
\wt w
=
\sqrt
{\fr{\wt\tau-\wt w(\wt\tau'-\wt w)}
{\wt\tau}}.
\ee
Let us introduce the function $j$ by means of the equality
\be
\tau=j\wt w^2.
\ee
We obtain
$
[\wt\tau-\wt w(\wt\tau'-\wt w)]/\wt\tau
=
[1-j- \wt w {\wt j}']/  j
$
\ses
and
\be
\lf(\fr1{\sin^2\varrho}   \D{\cos\varrho}W \rg)^2
\lf( \D W{\wt w}\wt w \rg)^2
=
\fr1j-1
-
\fr1j\D jW
 \D W{\wt w}\wt w.
\ee
Using (5.49) and (5.50) together with
\be
j=
 H^2
\fr1{1-W^2 }
\ee
\ses
(see (5.41)),
we can write
the equation
(5.94) in the form
$$
\lf(
\fr{1}{1-W^2} \D W{\wt w}\wt w
 \rg)^2
+
2W
\lf(
\fr{1}{1-W^2} \D W{\wt w}\wt w
 \rg)
 =
\fr {1}{H^2}
(1-W^2 )
-1,
 $$
which can conveniently be written as follows:
$$
\lf(
\fr{1}{1-W^2} \D W{\wt w}\wt w
+W \rg)^2
 =
\lf(\fr 1{H^2}-1\rg)
(1-W^2 ).
$$

It proves convenient to go over to the variable $W^2$:
\ses\\
$$
\lf(
\fr{1}{1-W^2} \D{ W^2}{\wt w}\wt w
+2W^2
 \rg)^2
 =
4
\lf(\fr 1{H^2}-1\rg)
(1-W^2 )
W^2.
$$
Since
\ses\\
$$
W^2=1-\fr{H^2\wt w^2}{\tau}
$$
(see (5.45)),
we get
$$
\lf(
-
\fr{\tau}{\wt w^2}
 \D{\fr{\wt w^2} {\tau} }{\wt w}\wt w
+2\lf(1-\fr{H^2\wt w^2}{\tau}\rg)
 \rg)^2
 =
 4(1-H^2)
\fr{\wt w^2}{\tau}
\lf(1-\fr{H^2\wt w^2}{\tau}\rg),
$$
\ses
or
$$
\lf(
\wt w\wt\tau'-2\tau
+2(\tau-H^2\wt w^2)
 \rg)^2
 =
 4(1-H^2)
\wt w^2
(\tau-H^2\wt w^2).
$$
Simplifying leaves us with the equation
$$
(\wt\tau' -2H^2\wt w )^2
 =
 4(1-H^2)
(\tau-H^2\wt w^2),
$$
which can readily be solved to yield (5.69).
Proposition 5.3 is valid.

\ses

The coefficients (5.72) show the properties
$$
u_kN^k{}_n=-\fr1h gqy^j \nabla_nb_j-
u_k  a^k{}_{nj}y^j,
\qquad
b_kN^k{}_n=\fr1h(1-h)y^j \nabla_nb_j
-
b_k  a^k{}_{nj}y^j
$$
(where $u_k=a_{kn}y^n$),
and
$$
d_nb\equiv
\D{b}{x^n}+b_kN^k{}_n=\fr1hy^j \nabla_nb_j,
\qquad
d_nq\equiv
\D{q}{x^n}+\fr1q v_kN^k{}_n=-\fr1{hq}
(b+gq)
 y^j \nabla_nb_j,
$$
together with
$$
d_n\lf(\fr qb\rg)
=
-\fr1{b^2qh}B
 y^j \nabla_nb_j,
 \quad
d_n B
=
-\fr g{qh}B
 y^j \nabla_nb_j,
   \quad
d_n\fr B{b^2}=
-\fr{2q+g b}
{b^3qh}B
 y^j \nabla_nb_j .
$$

\ses

With these formulas  it is possible to verify directly the validity of the
vanishing
\be
\cD_n K~:=\D K{x^n}+N^m{}_nl_m=0,
\ee
\ses
\be
\cD_n y_j~:=  \D {y_j}{x^n}+N^m{}_ng_{mj} - D^m{}_{nj}y_m=0,
\ee
and
\be
\cD_n g_{ij} ~:=\D {g_{ij}}{x^n}+2N^m{}_nC_{mji} - D^m{}_{nj}g_{mi}-
D^m{}_{ni}g_{mj}
=0,
\ee
where $D^m{}_{nj}=-\partial{N^m{}_{n}}/\partial {y^j}.$

\ses

\ses

The identity
\be
 l_hN^h{}_i=
-
K\fr{g\wt q}B(   y^h\nabla_i \wt b_h)
 -
l_t a^t{}_{ih}
y^h
 \ee
\ses
coming from (5.72) is useful to take into account when considering the vanishing (5.96).
The vanishing (5.97)  can be obtained directly by differentiating (5.96)
with respect to $y^j$.

Using (5.96) and (5.99), we can modify    the representation   (5.72)
by evaluating the sum
$$
N^m{}_i+l^m\D K{x^i}=  N^m{}_i-l^mN^h{}_il_h
=
\fr{g\wt q}By^m \wt s_i
+l^m   l_t a^t{}_{ih}  y^h
$$

\ses

$$
+
\fr1h
\fr1{\wt q}
Km^m
\wt s_i
-
\fr1h
\lf(b+\fr12g\wt q\rg)
\wt\beta^m_i
-
\wt b^m
\wt s_i
+
\wt b \nabla_i \wt b^m
 -
 a^m{}_{ih}
y^h.
$$
\ses
We insert here (5.74),  getting
$$
N^m{}_i=-l^m\D K{x^i}
+
\fr{g\wt q}By^m \wt s_i
+
\fr1h
\fr1{\wt q}
Km^m
\wt s_i
$$

\ses

$$
+\lf[\wt b-
\fr1h
\lf(b+\fr12g\wt q\rg)
\rg]
\nabla_i \wt b^m
+
\fr1h
\lf(b+\fr12g\wt q\rg)
            \fr1{ \wt q^2} \wt s_i
\lf(
y^m
-
\wt b \wt b^m
\rg)
-
\wt b^m
\wt s_i
-    h^m_t a^t{}_{ij}y^j.
$$

Let us introduce the  tensor
\be
\wt\eta^{kn}=a^{kn}-\wt b^k\wt b^n-\fr1{\wt q^2}\wt v^k\wt v^n, \qquad
\wt v^k=y^k-\wt b \wt b^k.
\ee

We come to
\ses\\
$$
N^m{}_i=-l^m\D K{x^i}
+
\fr{g\wt q}By^m \wt s_i
+
\fr1h
\fr1{\wt q}
Km^m
\wt s_i
$$

\ses
\ses

$$
+\lf[\wt b-
\fr1h
\lf(b+\fr12g\wt q\rg)
\rg]
\wt\eta^{mj}
\nabla_i \wt b_j
+\lf[\wt b-
\fr1h
\lf(b+\fr12g\wt q\rg)
\rg]
 \fr1{\wt q^2}\wt v^m
\wt s_i
$$

\ses

\ses

$$
+
\fr1h
\lf(b+\fr12g\wt q\rg)
            \fr1{ \wt q^2} \wt s_i
\wt v^m
-
\wt b^m
\wt s_i
-    h^m_t a^t{}_{ij}y^j.
$$

\ses

\ses

In this way,
with  the tensor
\be
{\cal H}^{mj}~:=g^{mj}-l^ml^j-m^mm^j,
\ee
\ses
  we arrive at   the representation
\be
N^m{}_i
=
-
l^m\D K{x^i}
+
\Biggl[\lf(\wt b-\fr1h\lf(b+\fr g2\wt q\rg)\rg)
{\cal H}^{mj}\fr{K^2}{ B}
+
\lf(\fr1{h\wt q} -\fr{b^2+\wt q^2}{\wt qB}\rg) Km^m
y^j
\Biggr]
 \nabla_i\wt b_j
-    h^m_t a^t{}_{ij}y^j,
\ee
where $h^m_t=\de^m_t-l^ml_t$
and $m^m$ is the vector (5.73).

The equality
${\cal H}^{mj} =(B/K^2)  \wt  \eta^{mj}$
holds.

In the dimension $N=2$ we would have
${\cal H}^{mj}=0.$

Regarding regularity of the global $y$-dependence,
it should be noted that
the ${\mathbf\cF\cF^{PD}_{g}}$-Finsleroid metric function $K$
given by the formulas
(5.76)-(5.80)
involves the scalar
$ \wt q=\sqrt{\wt r_{mn}y^my^n}$ with $\wt r_{mn}=a_{mn}-\wt b_m\wt b_n.
$
Since
the  1-form $\wt b$  is of the unit norm $||\wt b||=1$,
the scalar  $\wt q$ is zero when $y=b$ or $y=-b$,
that is, in the directions of the north pole or the south pole
of the Finsleroid.
The derivatives of $K$ may
involve the fraction $1/{\wt q}$
 which gives rise to the {\it pole singularities}  when $\wt q=0$.
   This  just happens in the
right-hand part of the representation
(5.102)
for the coefficients
$N^m{}_i$.

 Therefore, we may apply the coefficients on but
 the  $b$-slit  tangent bundle
$
{\cal T}_bM~:= TM \setminus0\setminus b\setminus -b
$
(obtained
by deleting out in $TM\setminus 0$
all the directions which point along, or oppose,
 the directions  given rise to by the  1-form $b$), on which
 the coefficients $N^m{}_i$, as well as the function $K$,
are  smooth  of the class
  $C^{\infty}$ regarding the $y$-dependence.

On  the punctured tangent bundle $TM\setminus0$,
the metric function $K$
 is smooth globally
 of the class $C^2$ and not of the class $C^3$
 regarding the $y$-dependence.

In the case  (5.70) the equation  (5.43) can readily be solved, yielding
\be
\rho =f,
\ee
where $f$ is the function which was  indicated in (5.83).
We obtain
\be
\sin\varrho=\fr{h\wt q}{{\sqrt B}}, \qquad
\cos\varrho=\fr{b+\fr12g\wt q}{{\sqrt B}}.
\ee
The representation (5.40) entails
\be
\mu=h^2,
\ee
so that from (5.26) we may conclude that $C_1=0$.
The transformation
(5.27)  reduces to
\be
  t^m
=
\lf[
h
(y^m- \wt b\wt b^m)
+
\lf(b+\fr12g\wt q\rg)
\wt b^m
\rg]
\fr{K^h}{\sqrt B}.
\ee

Thus we have

\ses

\ses

{\bf Proposition 5.4.}
{\it
In
 the ${\mathbf\cF\cF^{PD}_{g}}$-Finsleroid space
 the transformation} (5.106) {\it performs  the conformally automorphic
transformation.
When
 $h=const$   and    $c=const$,
  the coefficients} (2.36)
{\it can explicitly be given by means of
the representation}
(5.101)-(5.102).

\ses

\ses

In the remainder of the present section, we take $c=1$, that is,
$|| b||_{\text{Riemannian}}=1.
$
Using (5.73), we can transform (5.106) to the expansion
\be
  t^m =    (T_1 l^m + T_2  m^m)
  \fr{K^2}B
\fr{K^{h-1}}{\sqrt B}
\ee
with respect to the frame $\{l^m,m^m\}$,
where
\be
T_1=-(1-h)  q^2+B+  \fr12gq(b+g q),
\qquad
T_2=\lf((1-h)b+\fr12g q\rg) q.
\ee

\ses

The $t^m$ of (5.106) is equivalent to the $\zeta^m$ of (6.26) of [7]:
$
t^m\equiv \zeta^m.
$
The coefficients (5.102)  are  equivalent to (6.62) of [7].
Therefore, with the substitution
$
\zeta^m=t^m
$
 all the relations among curvature tensors
which were established in [7]
are applicable  to the approach developed in the present section,
including the following:
$$
\fr B{K^2} M_{nij}=
\Bigl((1-h)b  +\fr12g{ q}    \Bigr)
\fr1h
b_la_n{}^l{}_{ij}
-
\lf(
\fr g{2{ q}} v_n
+
(1-h)b_n
\rg)
\fr1h
 y^t
 b_la_t{}^l{}_{ij}
-
a_{tnij}
 y^t
$$
\ses
and
$$
\fr B{K^2}M^{nij} M_{nij}
=
\Biggl(\fr1{h}
\Bigl((1-h)b  +\fr g2  q    \Bigr)
b_ha^{nhij}
-
a_h{}^{nij}
 y^h
\Biggr)
\Biggl(\fr1{h}
\Bigl((1-h)b  +\fr g2  q    \Bigr)
b_la_n{}^l{}_{ij}
-
a_{tnij}
 y^t
\Biggr).
$$

\ses

\ses

If we take $\la$ from (2.32) and
 the coefficients
$N^k{}_i$ from (5.102),
and use the functions $t^m=t^m(x,y)$ specified by (5.106),
we obtain
the vanishing
$  d_i\la(x,y_1,y_2)=0,$
when $h=const.$
To verify the statement, it is worth
deriving  the equality
\be
\D{\la}{y^k_1}=
h^2\fr{
B_1v_{2k}+q^2_1b_kA_2
- b_1A_2
v_{1k}
-v_{12}
\lf(h^2v_{1k}+\lf(b_k+\fr12g\fr{1}{q_1}v_{1k}\rg)A_1
\rg)
}
 {B_1 \sqrt{B_1}\,\sqrt{B_2} },
\ee
together with the counterpart
\be
\D{\la}{y^k_2}=
h^2\fr{
B_2v_{1k}+q^2_2b_kA_1
- b_2A_1
v_{2k}
-v_{12}
\lf(h^2v_{2k}+\lf(b_k+\fr12g\fr{1}{q_2}v_{2k}\rg)A_2
\rg)
}
 {B_2 \sqrt{B_2}\,\sqrt{B_1} },
\ee
where
$
A_1=A(x,y_1),\,   A_2=A(x,y_2),\,
B_1=B(x,y_1),\,   B_2=B(x,y_2),\,
q_1=q(x,y_1),\,q_2=q(x,y_2),\,
b_1=b(x,y_1),\,   b_2=b(x,y_2),$
together with
$
v_{1i}= r_{in}(x)y^n_1$ and $  v_{2i}= r_{in}(x)y^n_2.$
Plugging these derivatives in
$  d_i\la(x,y_1,y_2)$ results in the claimed vanishing
$  d_i\la(x,y_1,y_2)=0$
after attentive couplepage reductions.

\ses

It will be noted that
$$
b^k\D{\la}{y^k_1}=  h^2\fr{  q^2_1A_2  -v_{12}   A_1   }
 {B_1 \sqrt{B_1}\,\sqrt{B_2} },  \qquad
 b^k\D{\la}{y^k_2}=  h^2\fr{  q^2_2A_1  -v_{12}   A_2   }
 {B_2 \sqrt{B_1}\,\sqrt{B_2} }.
$$

We have also
$$
\D{\la}g=-\fr12\lf(\fr {b_1q_1}{{B_1}}+\fr{b_2q_2}{{B_2}}\rg)  \la
+\fr{q_1A_2+q_2A_1-gv_{12}}
 {2 \sqrt{B_1}\,\sqrt{B_2} },
$$
or
$$
\D{\la}g=
\fr{1}
 { 2\sqrt{B_1}\,\sqrt{B_2} }
 \lf[
\fr{q_1^2 A_2}{B_1}\si_1
+
\fr{q_2^2A_1}{B_2}\si_2
-v_{12}
 \lf(
\fr{ A_1}{B_1}\si_1
+
\fr{ A_2}{B_2}\si_2
 \rg)
 \rg],
 $$
where
\be
\si_1=\fr g2A_1+h^2q_1\equiv q_1+\fr g2b_1,  \qquad   \si_2=\fr g2A_2+h^2q_2\equiv q_2+\fr g2b_2.
\ee
There arises the equality
\be
\D{\la}g=\fr1{2h^2}\lf[
\si_1b^k  \D{\la}{y^k_1}  +   \si_2b^k  \D{\la}{y^k_2}
\rg]
\equiv
\fr1{h^2}\lf[
z_1C_1^k  \D{\la}{y^k_1}  +   z_2C_2^k  \D{\la}{y^k_2}
\rg],
\ee
where
$$
z_1=   \fr{q_1K_1^2}{NgB_1}   \si_1,  \qquad
z_2=   \fr{q_2K_2^2}{NgB_2}   \si_2.
$$

\ses

\ses

\ses

\setcounter{sctn}{6} \setcounter{equation}{0}

\bc  { \bf 6. Conclusions}   \ec

\ses

\ses

In the two-dimensional approach, $N=2$,
the general representation  for the coefficients $N^m{}_i=N^m{}_i(x,y)$
entailing the  property of preservation of two-vector angle
can be indicated locally for arbitrary
sufficiently smooth Finsler metric function [8,9].
Such a general possibility can doubtfully  be meet in the dimensions $N\ge3$,
for in these dimensions the two-vector
is of complicated nature
 except for rare particular cases.
Such lucky cases are just proposed
 by the Finsler spaces
which are conformally automorphic to the  Riemannian spaces.
The respective two-vector angle is explicit,
namely is given by the simple formulas (1.7) and
(2.31)-(2.32).
Such Finsler spaces can be characterized by the constancy of the indicatrix curvature.
In each tangent space,
 the indicatrix curvature value
$
{\cal C}_{\text{Ind.}} = H^2
$
is obtained and the relevant  conformal multiplier is given by $p^2$ with
$
p=
(1/H)
F^{1-H}.
$
This $p$ is constructed from the Finsler metric function $F$.
The $H$ is
the degree of conformal automorphism.
In the case $H=1$ the Finsler space under consideration reduces to the Riemannian space proper.

In indicatrix-homogeneous case,
the required connection coefficients are presented by the pair
$\{N^j{}_i,\,D^j{}_{ik}\}$,
where
$D^j{}_{ik}=-\partial N^j{}_i/\partial y^k$.
The equality
$
N^j{}_i=-D^j{}_{ik}y^k
$
holds.

In the Riemannian geometry
the two-vector angle is
$\al^{\text{Riem}}_{\{x\}}(y_1,y_2)=a_{mn}(x)y^m_1y^n_2/S_1S_2,$
where $S_1=\sqrt{a_{mn}(x)y^m_1y^n_1}$ and $S_2=\sqrt{a_{mn}(x)y^m_2y^n_2}$.
Starting with the fundamental property of
 the
metrical  linear Riemannian  connection
 that the Riemannian
angle is preserving under the parallel displacements of the involved vectors, which in terms of our notation
can be written as
$$
d^{\text{Riem}}_i   \al^{\text{Riem}}_{\{x\}}(y_1,y_2)=0,
\qquad y_1,y_2\in T_xM,
$$
with
\ses
$$
d^{\text{Riem}}_i=\D{}{x^i}
+L^k{}_i{(x,y_1)} \fr{\partial}{\partial y_1^k}
+L^k{}_i{(x,y_2)} \fr{\partial}{\partial y_2^k},
$$
\ses\\
where
$
  L^k{}_i{(x,y_1)} =-a^k{}_{ij}(x)y_1^j,
$
\ses
$  L^k{}_i{(x,y_2)} =-a^k{}_{ij}(x)y_2^j,
$
and
   $a^k{}_{ij}$
     are the Riemannian    Christoffel symbols   fulfilling
the Riemannian Levi-Civita connection,
the important question can be set forth:
Can we  have  the similar vanishing  in the  Finsler   space?
It proves that  the respective extension of the Riemannian equation
$d^{\text{Riem}}_i   \al^{\text{Riem}}=0$
to the equation
$d_i   \al=0$
applicable to the Finsler space
under consideration
can  straightforwardly be solved
giving
the required
coefficients
$N^j{}_i$ indicated in
 (2.36).
They admit the remarkable alternative representation
$N^n{}_{i}=d^{\text{Riem}}_i   y^n$
(see (1.24)).
In this way we obtain the
connection $\{N^j{}_i,\,D^j{}_{ik}\}$
which is
metrical and simultaneously
 angle-preserving.
The key vanishing
$y_kN^k{}_{mnj}=0$
holds fine.

Remarkably, the Finsler connection presented by this pair
$\{N^j{}_i,\,D^j{}_{ik}\}$
is
the image of
the
metrical  linear Riemannian  connection
under conformally-automorphic transformations.
When going from the considered Finsler space to the underlined Riemannian space,
the covariant derivative behaves transitively
and the non-linear deformation which materializes the conformal automorphism
is parallel.
In particular,
the Riemannian vanishing
$
d^{\text{Riem}}_mS=0
$
\ses
 just entails the Finslerian counterpart
$
d_mF=0.
$

Also, the involved coefficients  $N^m{}_i$ fulfill  the  representation
$
N^k{}_{mnj}=  -   \cD_mC^k{}_{nj}
$
(see Proposition 3.2).
Just the same representation is valid in the two-dimensional Finsler spaces
(see (2.14)  in [8,9]).
Is the equation
$$
\Dd{N^k{}_{m}}{y^n}{y^j}=  -   \cD_mC^k{}_{nj}
$$
meaningful in other (in any?)
Finsler spaces to find the coefficients
$N^k{}_{m}$
required
to preserve the two-vector angle?
The question is addressed to readers.

The curvature tensor
$\Rho_k{}^n{}_{ij}$
has been explicated from commutators of arisen covariant derivatives
which is attractive to develop in future
the theory of curvature for the Finsler space
${\cal F}^N$.

For the  ${\cal FS}$-space specialized by (1.25)
 we have got at our disposal the simple
  example
  of the parallel deformation transformation, namely proposing by (5.27),
which  entails  the coefficients $N^m{}_i$
possessing the property of   angle preservation.
The coefficients  are  given
 explicitly
 by the representation (5.72)--(5.75), which admits the alternative form
(5.101)-(5.102).
The space proves to be  of the Finsleroid type,
with
 the Finsleroid characteristic parameter $g$
 manifesting the meaning:
$h=\sqrt{1-(g^2/4)}$ is  the homogeneity degree
(denoted above by $H$)
of the conformal automorphism.
The Finsleroid metric function $K$
when considered on the $b$-slit tangent bundle
$
{\cal T}_bM~:= TM \setminus0\setminus b\setminus -b
$
 is smooth of the class $C^{\infty}$
regarding the global $y$-dependence.
The same regularity property is valid for the coefficients
$N^m{}_i$
given by (5.102).

\ses

\ses

\setcounter{equation}{0}

\bc { \bf Appendix A: Proof of Proposition 2.1}  \ec

\ses

\ses

Let us verify the validity of  Proposition 2.1,
starting with
the conformal tensor
$$
u_{ij}=F^{2a}g_{ij}, \qquad a=a(x),
$$
and denoting
$
u_{ijk}=\partial u_{ij}/\partial{y^k}.
$
We get
$
u_{ijk}=2(a/F)F^{2a}l_kg_{ij}+2F^{2a}C_{ijk},
$
where
$C_{ijk}=(1/2)\partial g_{ij}/\partial y^k$.
Constructing the coefficients
$$
Z_{ijk}~:=\fr12(u_{kji}+u_{iki}-u_{ijk})
$$
leads to
$$
Z_{ijk}=\fr a{F}F^{2a}(l_ig_{kj}+l_jg_{ik}-l_kg_{ij})+F^{2a}C_{ijk}.
$$
Since the components $u^{ij}$ reciprocal to $u_{ij}$ are of the form
$
u^{ij}=F^{-2a}g^{ij},
$
the coefficients
$
Z^m{}_{ij}=u^{mh}Z_{ijh}
$
read merely
\ses\\
$$
Z^m{}_{ij}=\fr a{F}(l_i\de^m_j+l_j\de^m_i-l^mg_{ij})+C^m{}_{ij}.
$$

We obtain
$$
\D{Z^m{}_{ni}}{y^j}
=
-
\fr {a}{F^2}
l_j(l_n\de^m_i+l_i\de^m_n-l^mg_{ni})
+
\fr a{F^2}
(h_{ij}\de^m_n+h_{nj}\de^m_i-h^m_jg_{in}-2l^mFC_{inj})
+
\D{C^m{}_{ni}}{y^j}
$$
and
$$
\D{Z^m{}_{ni}}{y^j}-
\D{Z^m{}_{nj}}{y^i}
=
-
\fr {a}{F^2}
[l_n(l_j\de^m_i-l_i\de^m_j) -l^m(l_jg_{ni}-l_ig_{nj})]
$$

\ses

\ses

$$
+
\fr a{F^2}
[(h_{nj}\de^m_i-h_{ni}\de^m_j)  -(h^m_jg_{in}-h^m_ig_{jn})]
+
\D{C^m{}_{ni}}{y^j}
-
\D{C^m{}_{nj}}{y^i},
$$
\ses
so that
$$
\D{Z^m{}_{ni}}{y^j}-    \D{Z^m{}_{nj}}{y^i}
=
\fr {2a}{F^2}
(h_{nj}h^m_i-h_{ni}h^m_j)
+
\D{C^m{}_{ni}}{y^j}
-
\D{C^m{}_{nj}}{y^i}.
$$

\ses

Also,
\ses\\
$$
Z^h{}_{ni}    Z^m{}_{hj}
-
Z^h{}_{nj}    Z^m{}_{hi}
=
 \fr a{F}
\Bigl[
\fr a{F}[l_n(l_i\de^m_j-l_j\de^m_i) +l^m(l_jg_{ni}-l_ig_{nj})]
+(l_iC^m{}_{nj}-l_jC^m{}_{ni})
\Bigr]
$$

\ses

\ses

$$
-
\lf(\fr a{F}\rg)^2
(g_{in}\de^m_j-g_{jn}\de^m_i)
+
\fr a{F}
(l_jC^m{}_{in}-l_iC^m{}_{jn})
+
C^h{}_{ni}    C^m{}_{hj}
-
C^h{}_{nj}    C^m{}_{hi},
$$
\ses
or
$$
Z^h{}_{ni}    Z^m{}_{hj}
-
Z^h{}_{nj}    Z^m{}_{hi}
=
-
\lf(\fr a{F}\rg)^2
(h_{in}h^m_j-h_{jn}h^m_i)
+
C^h{}_{ni}    C^m{}_{hj}
-
C^h{}_{nj}    C^m{}_{hi}.
$$

\ses

The curvature tensor
$$
\wt R_n{}^m{}_{ij}~:=
\D{Z^m{}_{ni}}{y^j}-
\D{Z^m{}_{nj}}{y^i}
+Z^h{}_{ni}    Z^m{}_{hj}
-
Z^h{}_{nj}    Z^m{}_{hi}
$$
\ses
is found as follows:
$
F^2\wt R_n{}^m{}_{ij}
=
a(a+2)
(h_{nj}h^m_i-h_{ni}h^m_j)
+
S_n{}^m{}_{ij},
$
where
\ses\\
$$
S_n{}^m{}_{ij}=
\lf(
\D{C^m{}_{ni}}{y^j}
-
\D{C^m{}_{nj}}{y^i}
 +
C^h{}_{ni}    C^m{}_{hj}
-
C^h{}_{nj}    C^m{}_{hi}
\rg)
F^2.
$$
\ses\\
In term of the covariant components
$
\wt R_{nmij}
=
u_{mh}\wt R_n{}^h{}_{ij}
$ and
$
S_{nmij}
=
g_{mh}S_n{}^h{}_{ij},
$
\ses
we obtain
$$
F^2\wt R_{nmij}
=
S_{nmij}
+
a(a+2)
(h_{nj}h_{mi}-h_{ni}h_{mj}).
$$
\ses
Therefore,
if   $  \wt R_{nmij}=0,$
then
\ses\\
\be
S_{nmij}
=
C
(h_{nj}h_{mi}-h_{ni}h_{mj}), \qquad C=-a(a+2).
\ee
\ses\\
Since
$
{\cal C}_{\text{Ind.}}=1-C
$
(see Section 5.8 in [1]),
we get
$
{\cal C}_{\text{Ind.}}=H^2,
$
where $ H=a+1$.
The proposition is valid.

\ses

\ses

\setcounter{equation}{0}

\bc { \bf Appendix B: Proof of Proposition 2.2}  \ec

\ses

\ses

Let us verify the validity of  Proposition 2.2. From the equation
$$
\D{\la}{x^i}
+N^k{}_{1i}
\D{\la}{y^k_1}
+
N^k{}_{2i}
\D{\la}{y^k_2}
=0
$$
we want to find the tensors
\be
n^m_{1i}=t_{1k}^mN^k{}_{1i}, \qquad
n^m_{2i}=t_{2k}^mN^k{}_{2i}.
\ee
Using (2.33) and (2.34), we obtain
\ses\\
$$
  \fr{a_{mn,i}t_{1}^m  t_2^n}
{ S_1S_2 }
-
\fr12\la
\lf[
\fr1{S_1S_1 }
a_{mn,i}t_{1}^mt_1^n
    +
\fr1{S_2S_2}
a_{mn,i}t_{2}^mt_2^n
\rg]
$$

\ses

\ses

$$
+
\lf[
  \fr{a_{mn}  t_2^n}
{S_1S_2 }
 -      \fr{a_{mn} t_1^n}
{S_1S_1}\la
\rg]
\lf(n^m_{1i}+\D{t_{1}^m}{x^i}\rg)
+
\lf[
  \fr{a_{mn}  t_1^n}
{ S_2S_1 }
 -      \fr{a_{mn} t_2^n}
{S_2S_2}\la
\rg]
\lf(n^m_{2i}+\D{t_{2}^m}{x^i}\rg)
=0,
$$
\ses
which can be written in the concise form
$$
\fr1{S_1}
\lf[
  \fr{a_{mn}  t_2^n}
{ S_2 }
 -      \fr{a_{mn} t_1^n}
{S_1}\la
\rg]
\nu^m_{1i}
+
\fr1{S_2}
\lf[
  \fr{a_{mn}  t_1^n}
{S_1 }
 -      \fr{a_{mn} t_2^n}
{ S_2}\la
\rg]
\nu^m_{2i}
=0,
$$
where
$$
\nu^m_{1i}=n^m_{1i}+\D{t_{1}^m}{x^i}  +a^m{}_{ik}t_1^k , \qquad
\nu^m_{2i}=n^m_{2i}+\D{t_{2}^m}{x^i} +a^m{}_{ik}t_2^k ,
$$
and
$a^m{}_{ik}$ are the Riemannian Christoffel symbols (2.37).

In this way
we come to the equation
\be
\lf(S_1S_2
  a_{mn}  t_2^n
 -
 S_2S_2   a_{mn} t_1^n
\rg)
\nu^m_{1i}
+
\lf(
S_1S_2
  a_{mn}  t_1^n
 -
 S_1S_1    a_{mn} t_2^n
\rg)
\nu^m_{2i}
=0.
\ee

\ses

Use
$$
d_iF=\D F{x^i}+l_kN^k{}_i
=
\D F{x^i}+\fr1HF^{2(1-H)}t_mn^m{}_i,
$$
so that
\ses\\
$$
t_mn^m{}_i
=
HF^{2(H-1)}
\lf(d_iF-\D F{x^i}\rg).
$$
From $S^2=F^{2H}$ it follows that
$$
t_m\lf(\D{t^m}{x^i}   +a^m{}_{ik}t^k \rg)
=
H\fr1F F^{2H} \D F{x^i}
$$
($H=const$ is implied).
We obtain
\be
t_m \nu^m{}_i
=
HF^{2(H-1)}
d_iF,
\qquad
\nu^m{}_i=n^m{}_i+\D{t^m}{x^i}   +a^m{}_{ik}t^k,
\qquad
n^m_{i}=t_{k}^mN^k{}_{i} ,
\ee
where the equality
$t_ht^h_n=
HF^{2(H-1)}y_n
$
(see (2.16)) has been used.
When
$
d_iF=0,
$
we have
unambiguously
$t_m \nu^m{}_i=0$ and the equation (B.2) reduces to
\be
  a_{mn}  t_2^n
\nu^m_{1i}
+
  a_{mn}  t_1^n
\nu^m_{2i}
=0.
\ee

Thus we may conclude that when $H=const$
and
$
d_iF=0
$
is fulfilled, the started equation
$d_i\la=0$
is equivalent to the equation (B.4).

The case
$ \nu^m{}_i=0$
reads
 \be
n^m{}_{i}=-\lf(\D{t^m}{x^i} +  a^m{}_{ik}t^k \rg),
\ee
 which
is equivalent to (2.36).
The examined proposition is valid.

\ses

\ses

\setcounter{equation}{0}

\bc { \bf Appendix C: Validity of Proposition 3.1}  \ec

\ses

\ses

Let us consider  the term
$$
a_{sh}t^h_{ni}
         T^s_m
+
t_s
\lf(    \D{t^s_{ni}}{x^m}    +   a^s{}_{mh}t^h_{ni}   \rg)
=
\D{t_ht^h_{ni} }{x^m}
-
t^jt^h_{ni}
\D{a_{jh}}{x^m}
+
a_{sj}t^j_{ni}
a^s{}_{mh}t^h
+
t_s
   a^s{}_{mh}t^h_{ni}
$$

\ses

\ses

\ses

$$
=
\D{t_ht^h_{ni} }{x^m}
-
t^jt^h_{ni}
\D{a_{jh}}{x^m}
+
\fr12
t^j_{ni}
\lf(
\D{a_{jh}}{x^m}
+
\D{a_{jm}}{x^h}
-
\D{a_{mh}}{x^j}
\rg)
t^h+
\fr12
t^j
\lf(
\D{a_{jh}}{x^m}
+
\D{a_{jm}}{x^h}
-
\D{a_{mh}}{x^j}
\rg)
t^h_{ni}
$$

\ses

\ses

$$
=
\D{t_ht^h_{ni} }{x^m}.
$$
We can take  $t_ht^h_{ni} $
from  (2.18).
By doing so and introducing the notation $P=1-H$, we transform the representation (3.3) to
$$
y_kN^k{}_{mni}
+
2C_{lni}N^l{}_m
=
  g_{ki}
(
  y^k_{sl}t^l_n
        T^s_m
+  y^k_s
T^s_{n,m}
)
+
  g_{kn}
(
  y^k_{sl}t^l_i
         T^s_m
+  y^k_s
T^s_{i,m}
)
$$

\ses

\ses

$$
-
\fr{2 P}H
F^{-2H}
\lf[
    (g_{ni}-2Hl_nl_i)   t_s
+
(y_na_{sl}t^l_i+y_ia_{sl}t^l_n)
\rg]
         T^s_m
$$

\ses

\ses

$$
-
\lf(
2
\fr PH
F^{-2H}y_nt_s
+
\fr1HF^{2(1-H)}a_{sh}t^h_n
\rg)
T^s_{i,m}
   -
\lf(
2
\fr PH
F^{-2H}y_it_s
+
\fr1HF^{2(1-H)}a_{sh}t^h_i
\rg)
T^s_{n,m}
 $$

\ses

\ses

$$
   -
P
F^{2(1-H)}
   \D{F^{2(H-1)}(g_{ni}-2l_nl_i)}{x^m}
 $$
\ses
and take the tensor $g_{ki}$ from (2.9),
obtaining
$$
y_kN^k{}_{mni}
+
2C_{lni}N^l{}_m
=
  g_{ki}
  y^k_{sl}t^l_n
         T^s_m
+
  g_{kn}
  y^k_{sl}t^l_i
        T^s_m
$$

\ses

\ses

$$
-
\fr {2P}H
F^{-2H}
\lf[
    (g_{ni}-2Hl_nl_i)   t_s
+
(y_na_{sl}t^l_i+y_ia_{sl}t^l_n)
\rg]
         T^s_m
$$

\ses

\ses

$$
-
\lf(
2
\fr PH
F^{-2H}y_nt_s
-
\fr P{H^2}F^{2(1-H)}a_{sh}t^h_n
\rg)
T^s_{i,m}
   -
\lf(
2
\fr PH
F^{-2H}y_it_s
-
\fr P{H^2}
F^{2(1-H)}a_{sh}t^h_i
\rg)
T^s_{n,m}
$$

\ses

\ses

$$
   -
P
F^{2(1-H)}
   \D{\lf(\fr1{H^2}a_{su}t^u_nt^s_i-2F^{2(H-1)}l_nl_i\rg)}{x^m}.
 $$

\ses

After that, we take into account the formula (3.31) which specifies the object
$
 T^s_{n,m}.
 $
This yields
 $$
y_kN^k{}_{mni}
+
2C_{lni}N^l{}_m
=
  g_{ki}
  y^k_{sl}t^l_n
         T^s_m
+
  g_{kn}
  y^k_{sl}t^l_i
         T^s_m
$$

\ses

\ses

$$
-
\fr {2P}H
F^{-2H}
\lf[
    (g_{ni}-2Hl_nl_i)   t_s
+
(y_na_{sl}t^l_i+y_ia_{sl}t^l_n)
\rg]
         T^s_m
$$

\ses

\ses

$$
-
2
\fr PH
F^{-2H}t_s
\lf[
y_n
 \lf(\D{t^s_i}{x^m}    +   a^s{}_{mh}t^h_i   \rg)
+
y_i
 \lf(\D{t^s_n}{x^m}    +   a^s{}_{mh}t^h_n   \rg)
 \rg]
+
\fr 1{H^2}
F^{2(1-H)}a_{sl}
a^s{}_{mh}
(t^l_nt^h_i+t^l_it^h_n)
$$

\ses

\ses

$$
   -
 \fr P{H^2}
F^{2(1-H)}
 t^u_nt^s_i   \D{a_{su}}{x^m}
+
2
 \fr P{H^2}
F^{2(1-H)}
   \D{F^{-2H}t_ht_st^h_nt^s_i}{x^m}.
    $$
\ses\\
Here,
$ t_ht^h_n=
HF^{2(H-1)}y_n
$
(see (2.16)).

\ses

Noting that
$$
  g_{ki}
  y^k_{sl}t^l_n
=
p^2a_{uv}t^u_kt^v_i
\D{  y^k_{s}}{y^n}
=
-p^2a_{uv}t^v_it^u_{kn}y^k_s,
$$
 taking  $C_{lni}$  from (2.19),
and using the vanishing
$
Hy_ky^k_s
-
F^{2(1-H)}t_s
=0
$
(see (2.15)),
we perform simplifications and remain with
\ses\\
$$
y_kN^k{}_{mni}
=
4P
F^{-2H}
l_nl_i   t_s
         T^s_m
-
\fr{2 P}H
F^{-2H}
(y_na_{sl}t^l_i+y_ia_{sl}t^l_n)
         T^s_m
-
2
\fr PH
F^{-2H}t_s
a^s{}_{mh}
(y_nt^h_i+y_it^h_n)
$$

\ses

\ses

 $$
+
\fr 1{H^2}F^{2(1-H)}a_{sl}
a^s{}_{mh}(t^l_nt^h_i+t^l_it^h_n)
   -
 \fr P{H^2}
F^{2(1-H)}
 t^u_nt^s_i   \D{a_{su}}{x^m}
+
2
 \fr P{H^2}
F^{2(1-H)}
 t^h_nt^s_i  \D{F^{-2H}t_ht_s}{x^m}.
    $$

\ses

Finally, we apply  (3.4)
$$
\D F{x^m}=
\fr F{HS^2}t_s
T^s_m,
$$
\ses\\
take $T^s_m$ from (3.1), and notice the vanishing
$$
-
\fr{2 P}H
F^{-2H}
(y_nt^l_i+y_it^l_n)
    \D{t_l}{x^m}
+
2
 \fr P{H^2}
F^{2(1-H)}
 t^h_nt^s_i F^{-2H} \D{t_ht_s}{x^m}
=0.
 $$
\ses\\
We arrive at
\ses\\
$$
y_kN^k{}_{mni}
=
-
\fr{2 P}H
F^{-2H}
(y_nt^l_i+y_it^l_n)
\lf(       -t^s \D{a_{sl}}{x^m}  +a_{sl}a^s{}_{mh}t^h        \rg)
-
\fr{2 P}H
F^{-2H}t_s
a^s{}_{ml}
(y_nt^l_i+y_it^l_n)
$$

\ses

\ses

$$
+
\fr P{H^2}F^{2(1-H)}a_{su}
a^s{}_{mh}(t^u_nt^h_i+t^u_it^h_n)
   -
 \fr P{H^2}
F^{2(1-H)}
 t^u_nt^h_i   \D{a_{hu}}{x^m}
=0.
    $$

 We  have verified the validity of Proposition 3.1.

\ses

\ses

\setcounter{equation}{0}

\bc { \bf Appendix D: Verifying  Proposition 3.2}  \ec

\ses

\ses

With the convenient notation
$$
X^k{}_{mn}= y^k_{sl}t^l_n            a^s{}_{mh}t^h
+  y^k_s      a^s{}_{mh}t^h_n
-
y^k_s a^{su}  t^l_n \D{a_{ul}}{x^m}
$$
\ses
we can write (3.1) in the form
$$
N^k{}_{mn}
+
X^k{}_{mn}
=
 -  y^k_{sl}t^l_n
 \D{t^s}{x^m}
 -  y^k_s a^{su}   \D{a_{ul}t^l_n}{x^m}.
 $$
\ses\\
Differentiating this equality with respect to $y^j$
and using the notation
$X^k{}_{mnj}=\partial X^k{}_{mn}/\partial y^j$, we get
\ses\\
$$
N^k{}_{mnj}
+
X^k{}_{mnj}
=
 -  (y^k_{slj}t^l_n  +  y^k_{sl}t^l_{nj})
 \D{t^s}{x^m}
  -  y^k_{sl}t^l_n
 \D{t^s_j}{x^m}
   -
    y^k_{sl}t^l_j a^{su}   \D{a_{ul}t^l_n}{x^m}
 -  y^k_s a^{su}   \D{y^h_uZ^l_{nj}}{x^m},
 $$
\ses\\
where $Z^l{}_{nj}=a_{vl}t^v_ht^l_{nj}$
and the identity
$y^h_ut^v_h=\de^v_u$
has been taken into account.
Here, the equality
$
 y^k_s a^{su} y^h_u
 =
  p^2g^{kh}
$
should be used.

This method results in
\ses\\
$$
N^k{}_{mnj}
+
X^k{}_{mnj}
=
 -  (y^k_{slj}t^l_n  +  y^k_{sl}t^l_{nj})
 \D{t^s}{x^m}
  -  y^k_{sl}t^l_n
 \D{t^s_j}{x^m}
  -
   y^k_{sl}t^l_j a^{su}   \D{a_{ul}t^l_n}{x^m}
$$

\ses

\ses

$$
 -  y^k_s a^{su} a_{vl}t^v_ht^l_{nj}  \D{y^h_u}{x^m}
 - p^2g^{kh}
  \D{Z^l_{nj}}{x^m}.
  $$
Since
$$
C_{hnj}=
(1-H)\fr1{F}(l_jg_{hn}+l_ng_{hj}-l_hg_{nj})
+p^2 Z^l{}_{nj}
$$
(see (2.22)),
we can write
$$
N^k{}_{mnj}
+
X^k{}_{mnj}
=
 -
  \lf(\D{}{y^j}\D{y^k_l}{y^n}\rg)t^l_u
 y^u_s
 \D{t^s}{x^m}
  -  y^k_{sl}t^l_n
 \D{t^s_j}{x^m}
  -
y^k_{sl}t^l_j a^{su}   \D{a_{ul}t^l_n}{x^m}
$$

\ses

\ses

$$
 -  y^k_s a^{su} a_{vl}t^v_ht^l_{nj}  \D{y^h_u}{x^m}
+ g^{kh}
  a_{vl}t^v_ht^l_{nj}
  \D{p^2}{x^m}
-g^{kh}
\D{C_{hnj}}{x^m}
  $$

\ses

\ses

$$
-
(1-H)g^{kh}\fr1{F^2}(l_jg_{hn}+l_ng_{hj}-l_hg_{nj})\D  F{x^m}
$$

\ses

\ses

$$
+
(1-H)g^{kh}\fr1{F}
\lf[
\lf(\D{l_j}{x^m}g_{hn}+\D{l_n}{x^m}g_{hj}
-\D{l_h}{x^m}g_{nj}\rg)
+
\lf(l_j\D{g_{hn}}{x^m}+l_n\D{g_{hj}}{x^m}
-l_h\D{g_{nj}}{x^m}\rg)
\rg].
$$

Considering the relation
$$
 \lf(\D{}{y^j}\D{y^k_l}{y^n}\rg)t^l_u
 =
 \D{}{y^j}\lf(t^l_u\D{y^k_l}{y^n}\rg)
-
 \D{y^k_l}{y^n}t^l_{uj}
=
-\D{}{y^j}\lf(y^k_lt^l_{un}\rg)
-
 \D{y^k_l}{y^n}t^l_{uj}
   $$
and noting that
$
y^k_l=g^{kh}p^2a_{lv}t^v_h,
$
we obtain the useful equality
\ses\\
$$
 \lf(\D{}{y^j}\D{y^k_l}{y^n}\rg)t^l_u
=
-\D{}{y^j}\lf( g^{kh}p^2a_{lv}t^v_h t^l_{un}\rg)
-
 \D{y^k_l}{y^n}t^l_{uj}.
$$

Along this way we come to
\ses\\
$$
N^k{}_{mnj}
+
X^k{}_{mnj}
+
g^{kh}
\D{C_{hnj}}{y^u}
   y^u_s
a^s{}_{mh}t^h
$$

\ses

\ses

\ses

$$
=
\lf(
\D{}{y^j}\lf( g^{kh}p^2a_{lv}t^v_h t^l_{un}\rg)
+
 \D{y^k_l}{y^n}t^l_{uj}
\rg)
 y^u_s
 \D{t^s}{x^m}
  -  y^k_{sl}t^l_n
 \D{t^s_j}{x^m}
  -
y^k_{sl}t^l_j a^{su}   \D{a_{ul}t^l_n}{x^m}
$$

\ses

\ses

$$
 +
 p^2g^{kh} a_{sl}t^l_{nj}  \D{t^s_h}{x^m}
+ g^{kh}
  a_{vl}t^v_ht^l_{nj}
  \D{p^2}{x^m}
    $$

\ses

\ses

$$
-
(1-H)g^{kh}
\lf[
\fr1{F^2}(l_jg_{hn}+l_ng_{hj}-l_hg_{nj})\D  F{x^m}
-
\fr1{F}\lf(\D{l_j}{x^m}g_{hn}+\D{l_n}{x^m}g_{hj}
-\D{l_h}{x^m}g_{nj}\rg)
\rg]
  $$

\ses

\ses

$$
+
(1-H)g^{kh}\fr1{F}\lf(l_j\D{g_{hn}}{x^m}+l_n\D{g_{hj}}{x^m}
-l_h\D{g_{nj}}{x^m}\rg)
-
g^{kh}
\D{C_{hnj}}{y^u}
   y^u_s
         \D{t^s}{x^m}
-g^{kh}
d_mC_{hnj},
$$
\ses
where
$$
d_mC_{hnj}=
\D{C_{hnj}}{x^m}   +   N^u{}_m\D{C_{hnj}}{y^u}.
$$

Reducing similar terms yields
\ses\\
\be
N^k{}_{mnj}=
-g^{kh}
d_mC_{hnj}+J^k{}_{mnj},
\ee
\ses\\
where
$$
J^k{}_{mnj}=
-
X^k{}_{mnj}
-
g^{kh}
\D{C_{hnj}}{y^u}
   y^u_s
a^s{}_{mh}t^h
-
y^k_{sl}t^l_j a^{su}  t^l_n \D{a_{ul}}{x^m}
 $$

\ses

\ses

\ses

$$
+
\lf(
t^l_{un}
\D{}{y^j}\lf( g^{kh}p^2a_{lv}t^v_h \rg)
+
 \D{y^k_l}{y^n}t^l_{uj}
\rg)
 y^u_s
 \D{t^s}{x^m}
  -  y^k_{sl}t^l_n
 \D{t^s_j}{x^m}
  -
y^k_{sl}t^l_j    \D{t^s_n}{x^m}
$$

\ses

\ses

$$
 +
 p^2g^{kh} a_{sl}t^l_{nj}  \D{t^s_h}{x^m}
+ g^{kh}
  a_{vl}t^v_ht^l_{nj}
  \D{p^2}{x^m}
-
(1-H)g^{kh}\fr1{F^2}(l_jg_{hn}+l_ng_{hj}-l_hg_{nj})\D  F{x^m}
$$

\ses

\ses

$$
+
(1-H)g^{kh}\fr1{F}
\lf[
\lf(\D{l_j}{x^m}g_{hn}+\D{l_n}{x^m}g_{hj}
-\D{l_h}{x^m}g_{nj}\rg)
+
\lf(l_j\D{g_{hn}}{x^m}+l_n\D{g_{hj}}{x^m}
-l_h\D{g_{nj}}{x^m}\rg)
\rg]
  $$

\ses

\ses

\ses

\ses

$$
+
(1-H)g^{kh}\fr1{F^2}
\lf[
(l_jg_{hn}+l_ng_{hj}-l_hg_{nj})l_u
   y^u_s
 \D{t^s}{x^m}
-
\lf(h_{ju}g_{hn}+h_{nu}g_{hj}
-h_{hu}g_{nj}\rg)
   y^u_s
 \D{t^s}{x^m}
\rg]
   $$

\ses

\ses

\be
-
2(1-H)g^{kh}\fr1{F}\lf(l_jC_{hnu}+l_nC_{hju}   -l_hC_{nju}\rg)
   y^u_s
 \D{t^s}{x^m}
-
g^{kh}
t^l_{jn}
\D{}{y^u}\lf( p^2a_{lv}t^v_h \rg)
   y^u_s
         \D{t^s}{x^m}.
\ee

We also find the contraction
\ses\\
$$
N^r{}_{mh}C_{rnj}=
-
\lf[
  y^r_{sl}t^l_h
\lf(         \D{t^s}{x^m}  +a^s{}_{mh}t^h        \rg)
+
  y^r_s
\lf(    \D{t^s_h}{x^m}    +   a^s{}_{mv}t^v_h   \rg)
\rg]
C_{rnj}
$$
\ses
(see (3.1)),  getting
$$
N^r{}_{mh}C_{rnj}=
-
C_{rnj}
  y^r_{sl}t^l_h
  \D{t^s}{x^m}
-
\lf[
(1-H)\fr1{F}(l_jg_{rn}+l_ng_{rj}-l_rg_{nj})
+p^2 t^v_rt^l_{nj} a_{vl}
\rg]
  y^r_s
    \D{t^s_h}{x^m}
    $$

\ses

\ses

\be
-
\lf(
  y^r_{sl}t^l_h
a^s{}_{mh}t^h
+
  y^r_s
   a^s{}_{mv}t^v_h
\rg)
C_{rnj}.
\ee

\ses

The indicated relations are sufficient to obtain
the representation
$$
N^k{}_{mnj}=
-
\cD_mC^k{}_{nj}
\equiv
-d_mC^k{}_{nj} + N^k{}_{mt}C^t{}_{nj}
- N^t{}_{mn}C^k{}_{tj}
- N^t{}_{mj}C^k{}_{nt}
$$
after performing required substitutions.

\ses

\ses

\setcounter{equation}{0}

\bc { \bf Appendix E: Validity of Proposition 5.1}  \ec

\ses

\ses

With (5.39), we get the expression
\ses\\
$$
\fr{\mu}{H^2}
a_{mn} t^m_k   t^n_h  =
\lf[1-\sin^2\varrho
+
\fr 1{4 \breve k^2}
(1+\breve k^2  )^2
  \sin^2\varrho
\rg]
(\varrho')^2
\fr1{b^2}w^2 e_ke_h
F^{2H}
 $$

\ses

\ses

\ses

$$
+
  \sin^2\varrho F^{2H}
\lf(a_{hk}-\fr1{c^2}  b_kb_h
-
\fr1{\wt w^2}
[w^2e_k(w^2e_h + \wt w^2b_h)
 +
 \wt w^2b_k  (w^2e_h + \wt w^2b_h)
]
\rg)
\lf(\fr 1{b\wt w}\rg)^2
$$

\ses

\ses

\ses

$$
+
\sqrt{\mu}
\fr 1Fl_k
\lf[
\sqrt{\mu}
\fr 1Fl_h
+
\fr 1{2 \breve k}
(1-\breve k^2  )
 \sin\varrho\,
\varrho'
\fr wbe_h
\rg]
F^{2H}
$$

\ses

\ses

$$
+
\fr 1{2 \breve k}
(1-\breve k^2  )
 \sin\varrho\,
\varrho'
\fr wbe_k
\lf[
\sqrt{\mu}
\fr 1Fl_h
+
\fr 1{2 \breve k}
(1-\breve k^2  )
 \sin\varrho\,
\varrho'
\fr wbe_h
\rg]
F^{2H}
$$

\ses

\ses

\ses

 $$
+
 \sin\varrho
 \lf[\cos\varrho
-
\fr{1+\breve k^2  }{4 \breve k^2}
[
1-\breve k^2 +(1+\breve k^2  )   \cos\varrho]
\rg]
e_k
\lf[
\sqrt{\mu}
\fr 1Fl_h
+
\fr{1-\breve k^2  }{2 \breve k}
 \sin\varrho\,
\varrho'
\fr wbe_h
\rg]
\varrho'
\fr wb
\fr{H}{\sqrt{\mu} }
F^{2H}
$$

\ses

\ses

 $$
+
 \sin\varrho
 \lf[\cos\varrho
-
\fr{1+\breve k^2  }{4 \breve k^2}
[
1-\breve k^2 +(1+\breve k^2  )   \cos\varrho]
\rg]
e_h
\lf[
\sqrt{\mu}
\fr 1Fl_k
+
\fr{1-\breve k^2  }{2 \breve k}
 \sin\varrho\,
\varrho'
\fr wbe_k
\rg]
\varrho'
\fr wb
\fr{H}{\sqrt{\mu} }
F^{2H}
$$
which  can readily be simplified to read
\ses\\
$$
\fr{\mu}{H^2}
a_{mn} t^m_k   t^n_h  =
\lf[1
+
\fr 1{2 \breve k^2}
(1-\breve k^2  )^2
  \sin^2\varrho
\rg]
(\varrho')^2
\fr1{b^2}w^2 e_ke_h
F^{2H}
 $$

\ses

\ses

\ses

$$
+
  \sin^2\varrho F^{2H}
\lf(
a_{kh}-b_kb_h
-
w^2(e_h+b_h)(e_k+b_k)
+\lf(w^2-\fr{w^4}{\wt w^2}\rg)
e_ke_h
\rg)
\lf(\fr 1{b\wt w}\rg)^2
$$

\ses

\ses

\ses

$$
+
\mu
\fr 1{F^2}l_kl_h
+
\sqrt{\mu}
\fr1F
\fr 1{2 \breve k}
(1-\breve k^2  )
 \sin\varrho\,
\varrho'
\fr wb(l_ke_h+l_he_k)
F^{2H}
$$

\ses

\ses

\ses

 $$
-
 \fr{1-\breve k^2  }{4 \breve k^2}
 \sin\varrho
 \lf[
1+\breve k^2 +(1-\breve k^2  )   \cos\varrho
\rg]
\lf[
\fr {\sqrt{\mu}}F(e_kl_h+e_hl_k)
+
\fr{1-\breve k^2  }{ \breve k}
 \sin\varrho\,
\varrho'
\fr wb
e_ke_h
\rg]
\fr wb
\fr{\varrho'H}{\sqrt{\mu} }
F^{2H}
$$

\ses

\ses

\ses

\ses

$$
=
(\varrho')^2
\fr{w^2}{b^2} e_ke_h
F^{2H}
+
\mu
\fr 1{F^2}l_kl_h
 $$

\ses

\ses

\ses

$$
+
  \sin^2\varrho F^{2H}
\lf(
a_{kh}-b_kb_h
-
w^2(e_h+b_h)(e_k+b_k)
+\lf(w^2-\fr{w^4}{\wt w^2}\rg)
e_ke_h
\rg)
\lf(\fr 1{b\wt w}\rg)^2.
$$

\ses

We use here the expansion
$$
g_{kh}=
l_kl_h
+
\fr {F^{2}}{b^2\tau}
\lf(
a_{kh}-b_kb_h
-
w^2(e_h+b_h)(e_k+b_k)
+
\fr{\tau- w(\tau'- w)}   {\tau}
w^2e_he_k
\rg)
$$
of the involved Finslerian metric tensor
and take $\mu$ from (5.40) and
$
\wt w^2/\tau
$
from (5.41).

We obtain
\ses
$$
\fr{1}{\wt w^2}\tau\sin^2\varrho
\fr{1}{H^2}
a_{mn} t^m_k   t^n_h  =
(\varrho')^2
\fr{w^2}{b^2} e_ke_h
F^{2H}
+
\fr1{\wt w^2}
\tau  \sin^2\varrho F^{2(H-1)}
g_{kh} $$

\ses

\ses

\ses

\be
+
  \sin^2\varrho F^{2H}
\lf(
-
 \fr{\tau- w(\tau'- w)}   {\tau}
w^2e_he_k
+\lf(w^2-\fr{w^4}{\wt w^2}\rg)
e_ke_h
\rg)
\lf(\fr 1{b\wt w}\rg)^2.
\ee
This representation is obviously equivalent to (5.42).

\ses

\ses

\ses

\ses

\def\bibit[#1]#2\par{\rm\noindent\parskip1pt
                     \parbox[t]{.05\textwidth}{\mbox{}\hfill[#1]}\hfill
                     \parbox[t]{.925\textwidth}{\baselineskip11pt#2}\par}

\bc {\bf  References}
\ec

\ses

\bibit[1] H. Rund, \it The Differential Geometry of Finsler
 Spaces, \rm Springer, Berlin 1959.

\bibit[2] D. Bao, S. S. Chern, and Z. Shen, {\it  An
Introduction to Riemann-Finsler Geometry,}  Springer, N.Y., Berlin 2000.

\bibit[3] L. Kozma and L. Tam{\' a}ssy,
Finsler geometry without line elements faced to applications,
{   \it Rep. Math. Phys.} {\bf 51} (2003), 233--250.

\bibit[4] L. Tam{\' a}ssy,
 Metrical almost linear connections in $TM$ for Randers spaces,
{\it Bull. Soc. Sci. Lett. Lodz Ser. Rech. Deform } {\bf 51} (2006), 147-152.

\bibit[5] Z. L. Szab{\' o},  All regular Landsberg metrics are Berwald,
{\it Ann Glob Anal Geom  } {\bf 34} (2008), 381-386.

\bibit[6] L. Tam{\' a}ssy,  Angle in Minkowski and Finsler spaces,
{\it Bull. Soc. Sci. Lett. Lodz Ser. Rech. Deform } {\bf 49} (2006), 7-14.

\bibit[7] G. S. Asanov,
 Finsleroid  gives rise to the angle-preserving  connection,
 {\it  arXiv:} 0910.0935 [math.DG],  (2009).

\bibit[8] G. S. Asanov,
 Finslerian angle-preserving  connection in two-dimensional
case.  Regular realization,
 {\it  arXiv:} 0909.1641v1 [math.DG] (2009).

\bibit[9] G. S. Asanov,
 Finsler space connected by angle  in two dimensions. Regular case,
{\it Publ. Math. Debrecen } {\bf 77/1-2} (2010), 245--259.

\end {document}